%% file: arxiv_thesis_full.tex
\documentclass[11pt,
a4paper
,notitlepage
,abstract
,english
]
{scrartcl}

\usepackage{geometry}
\usepackage[T1]{ fontenc}
\usepackage[utf8]{inputenc}
\usepackage{wrapfig}
\usepackage[english]{babel}
\usepackage{amsmath}
\usepackage{amsfonts}
\usepackage{amssymb}
\usepackage{enumerate}
\usepackage[thmmarks]{ntheorem}
\usepackage{array}
\usepackage[final]{graphicx}
\usepackage{caption}
\usepackage{subcaption}
\usepackage{setspace}
\usepackage{arydshln}
\usepackage{hyperref}
\usepackage{bbold}
\usepackage{microtype}
\usepackage{lscape}

\newcommand{\CC}{\mathbb{C} }

\newcommand{\RR}{\mathbb{R} }

\newcommand{\PP}{\mathbb{P} }
\newcommand{\ZZ}{\mathbb{Z} }

\newcommand{\betrag}[1]{\left\lvert #1 \right\rvert}
\newcommand{\norm}[1]{\left\lVert #1 \right\rVert}

\newcommand{\anfzeichen}[1]{``#1''}
\newcommand{\diff}{\mathop{}\!\mathrm{d}}
\newcommand{\Diff}{\mathop{}\!\mathrm{D}}
\newcommand{\ie}{i.~e.\ }

\newcommand\numberthis{\addtocounter{equation}{1}\tag{\theequation}}

\theoremstyle{break}
\theoremseparator{:}
\newtheorem{theorem}{Theorem}[section]
\newtheorem{lemma}[theorem]{Lemma}
\newtheorem{proposition}[theorem]{Proposition}
\newtheorem{corollary}[theorem]{Corollary}

\theorembodyfont{\upshape}
\newtheorem{definition}[theorem]{Definition}
\newtheorem{remark}[theorem]{Remark}

\theoremstyle{nonumberplain}
\theoremheaderfont{\itshape}
\theoremsymbol{\ensuremath{\square}}
\newtheorem{proof}{Proof}

\title{Second species orbits of negative action and contact forms in the circular restricted three-body problem}
\author{Robert Nicholls}

\begin{document}
\maketitle

\begin{abstract}
\input{abstract.tex}
\end{abstract}

\tableofcontents

\section{Introduction} \label{chapter introduction}
\input{introduction}

\section{Energy hypersurfaces} \label{chapter hypersurfaces}
\input{hypersurfaces}

\section{Generating orbits} \label{chapter generating orbits}
\input{generating_orbits}

\section{Action of generating orbits} \label{chapter action of generating orbits}
\input{action}

\section{Proof of main theorem} \label{chapter proof of main theorem}
\input{main_proof}

\section{Numerical results} \label{chapter numerical}
\input{numerical}

\addcontentsline{toc}{section}{References}

\bibliography{paper}
\bibliographystyle{alpha}

\end{document}

%% file: abstract.tex
We show in this work that the restricted three-body problem is in general not of contact type and zero is the energy value where the contact property breaks down.
More explicitly, sequences of generating orbits with increasingly negative action and energies between~$-\sqrt{2}$ and zero are constructed.
Using results from~\cite{BolotinMacKay}, it is shown that these generating orbits extend to periodic solutions of the restricted three-body problem for small mass ratios and the action remains within a small neighbourhood.
These orbits obstruct the existence of contact structures for energy level sets~$\Sigma_c$ of the mentioned values and small mass ratios of the spatial problem.
In the planar case the constructed orbits are noncontractible even in the Moser-regularised energy hypersurface~$\overline{\Sigma}_c$.
Here, the constructed orbits still obstruct the existence of contact structures in certain relative de Rham classes of~$\overline{\Sigma}_c$ to the Liouville 1-form.
These results are optimal in the sense that for energies above zero the level sets are again contact for all mass ratios.
Numerical results are additionally given to visualise the computations and give evidence for the existence of these orbits for higher mass ratios.

%% file: introduction.tex

In order to use modern mathematical methods, often a contact structure is required.
A first step in connecting the restricted three-body problem to Reeb dynamics was done by Albers, Frauenfelder, van Koert and Parternain in~\cite{AlbersFrauenfelderVanKoertPaternain}, showing that the bounded components of the Moser-regularised energy hypersurface of the planar restricted three-body problem is of contact type below the first critical energy level and also slightly above this value.
The same result for the spatial problem was proven by Cho, Jung and Kim in~\cite{ChoJungKim_spatial}.
Although the spatial case has more physical relevance, for example in space mission design, some results on contact manifolds currently only work for three-dimensional manifolds, e.g.\ energy level sets inside a four-dimensional phase space of a two-dimensional configuration space.
For energies above zero the restricted three-body problem again admits a contact structure as the canonical Liouville 1-form becomes a contact form.
This can be checked easily by checking that the corresponding Liouville vector field is transverse to the energy hypersurfaces.

So the question remains:
What happens for energies in between these two values?
This present work focuses on the region between the highest critical energy value and zero, and the main statements are the following:

\begin{theorem}
	The spatial restricted three-body problem is in general not of contact type for energies between $-\sqrt{2}$ and zero.
\end{theorem}

\begin{theorem}
	If the planar restricted three-body problem is of contact type between~$-\sqrt{2}$ and zero, then the de Rham class of the contact form minus the Liouville 1-form must become infinitely bad for small mass ratios~$\mu$.
\end{theorem}

A more detailed version of the statements can be found in theorem~\ref{thm main theorem final}.
In order to explain the result more explicitly, let~$(T^*N, \diff \lambda)$ be the exact symplectic manifold of a cotangent bundle together with the exterior derivative of the canonical Liouville 1-form~$\lambda$.
To define the restricted three-body problem, let~$\mu$ bet the mass ratio of the two primaries~$M_1$ and~$M_2$.
The Hamiltonian is given by
\begin{align}
	H_\mu (q,p)
	=
	\frac{1}{2}\left( (p_1 + q_2)^2 + (p_2 - q_1)^2 + p_3^2 \right)- \frac{1 - \mu}{\norm{q - M_1}} - \frac{\mu}{\norm{q - M_2}} - \frac{q_1^2 + q_2^2}{2},
\end{align}
where~$q$ are the position coordinates in~$\RR^3 \setminus \{ M_1, M_2 \}$ and~$p$ are the momentum coordinates in~$T^*_q (\RR^3 \setminus \{ M_1, M_2 \})$.
The planar case is recovered by setting $q_3 = p_3 = 0$.
Dynamics are given by Hamilton's equation of motion and by the conservation of energy, solutions to this Hamiltonian system with energy $H_\mu = c$ stay in their own energy level set~$\Sigma_{\mu, c}$ for all time.
For all energies the hypersurface~$\Sigma_{\mu, c}$ can be regularised at collision with the primaries using the Moser regularisation.
We will denote the regularised energy hypersurface by~$\overline{\Sigma}_{\mu, c}$.
The question is whether one can find a contact structure~$\alpha$ on~$\overline{\Sigma}_{\mu, c}$ that is compatible with the symplectic structure~$\omega$ on~$T^* (\RR^n \setminus \{ M_1, M_2 \})$ in the sense that~$\diff \alpha = \omega$ and orientations are preserved.
If that is the case, the Reeb flow is a positive reparametrisation of the Hamiltonian flow and we can compute the integral of~$\alpha$ along a contractible orbit~$\gamma$ by
\begin{align*}
	0
	<
	\int \gamma^* \alpha
	=
	\int \gamma^* \lambda - \sum \left( r_i \int \gamma^* \beta_i \right),
\end{align*}
where~$\lambda$ is the canonical Liouville 1-form on the cotangent bundle, $\beta_i$ are de Rham generators and~$r_i$ are the corresponding coefficients of the de Rham class of~$\lambda - \alpha$ in~$\overline{\Sigma}_c$.
The integral~$\int \gamma^* \lambda$ is called the action of the orbit.
If one now has a contractible periodic orbit with negative action there can not exist a contact structure.
The action of noncontractible orbits, on the other hand, only obstructs the existence of contact structures in certain relative de Rham classes.

Over the course of this work we will first compute de Rham generators of the regularised energy hypersurfaces in chapter~\ref{chapter hypersurfaces}.
This is done by analysing the Moser regularisation in the setting of the theorem of Seifert-van Kampen and then applying these relations locally for the restricted three-body problem.

Then we will recall the notion of generating orbits from~\cite{Henon1} in chapter~\ref{chapter generating orbits}.
In short these generating orbits are limit orbits as~$\mu \to 0$.
They turn out to be either orbits of the limit Hamiltonian system, called the \emph{rotating Kepler problem}, or the Kepler problem in rotating coordinates~$(q,p)$, with its Hamiltonian
\begin{align} \label{eq rotating Kepler Hamiltonian}
	\begin{split}
		H_0 (q,p)
		&=
		H_\mathrm{fix}(q,p) + L(q,p)
		\\
		&=
		\frac{\norm{p}^2}{2} + q_1 p_2 - q_2 p_1 - \frac{1}{\norm{q}}
	\end{split}
\end{align}
or pieces of such solutions, glued together at collision with~$M_2$.

The advantage of working with these generating orbits is that one can now use the Kepler problem
\begin{align} \label{eq Kepler Hamiltonian}
	H_\mathrm{fix}(Q,P) = \frac{1}{2}\norm{P}^2 - \frac{1}{\norm{Q}},
\end{align}
where we know that solutions are ellipses with focus at the origin with period
\begin{align} \label{eq Kepler's third law}
	T &= 2 \pi \sqrt{a^3},
\end{align}
where~$a$ is the semi-major axis of the ellipse.
Furthermore, one can recover the Kepler energy~$H_\mathrm{fix}$ and the angular momentum~$L$ by
\begin{align} \label{eq Kepler energy in a}
	H_\mathrm{fix} &= - \frac{1}{2a}
	\quad \text{and}
	\\
	\label{eq angular momentum in a and epsilon'}
	L &= - \epsilon' \sqrt{a(1-\epsilon^2)},
\end{align}
where~$\epsilon$ is the eccentricity and~$\epsilon'$ the direction of rotation defined by
\begin{align} \label{eq definition epsilon prime}
	\epsilon' :=
	\begin{cases}
		+1 \quad \text{for anti-clockwise motion and}
		\\
		-1 \quad \text{for clockwise motion.} 
	\end{cases}
\end{align}
We call the case of~$\epsilon' = +1$ \emph{direct} or \emph{prograde motion} and $\epsilon' = -1$ \emph{retrograde}.

This is used in chapter~\ref{chapter action of generating orbits} to compute the action of generating orbits and ultimately find sequences of orbits with action tending towards negative values and energy tending to values between~$-\sqrt{2}$ and zero.
Another main ingredient for this computation is the Levi-Civita regularisation given by the map
\begin{align}\label{eq Levi-Civita transformation}
	\begin{split}
		l \colon \CC \setminus \{0\} &\to \CC \setminus \{0\}
		\\
		X &\mapsto X^2
	\end{split}
\end{align}
in complex coordinates~$X \in \CC \cong \RR^2$.
This maps solutions for the harmonic oscillator with spring constant $8 H_\mathrm{fix}$ onto Kepler solutions as a double cover with reparametrisation
\begin{align} \label{eq Levi-Civita time transformation}
	\diff t = 4\betrag{X}^2 \diff s
\end{align}
between the usual time~$t$ and the regularised time~$s$.

The corresponding solutions to the attractive harmonic oscillator are ellipses with centre at the origin and frequency
\begin{align} \label{eq regularised angular frequency}
	\varpi
	=
	\sqrt{-8 H_\mathrm{fix}}
\end{align}
for negative energies.
After rotation and time-shift as in~\cite{Celetti} we have solutions
\begin{align} \label{eq parametrization of regularised Kepler orbit}
	\begin{split}
		X_1(s) = \alpha \cos ( \varpi s)
		\qquad
		X_2(s) = \beta \sin ( \varpi s),
	\end{split}
\end{align}
where we can express the coefficients by
\begin{align}
	\label{eq alpha in a and epsilon}
	\alpha
	&=
	\sqrt{a(1+\epsilon)}
	\\
	\label{eq beta in a and epsilon}
	\beta
	&=
	\epsilon' \sqrt{a(1-\epsilon)}.
\end{align}
As a last ingredient we need the elapsed time of a Keplerian arc from collision to collision, which can be computed using Lambert's theorem and the free-fall time
\begin{align} \label{eq free-fall time}
	t_{Q_0}(Q_1)
	&=
	\sqrt{\frac{Q_0^3}{2}} \left( \sqrt{\frac{Q_1}{Q_0} \left( 1 - \frac{Q_1}{Q_0} \right) } + \arccos \left( \sqrt{\frac{Q_1}{Q_0}} \right) \right)
\end{align}
from height~$Q_0$ down to~$Q_1$.

Putting all the statements together, we can prove the main theorem in chapter~\ref{chapter proof of main theorem}.
We also add some numerical computations as visualisation and as evidence of how far these orbits survive in terms of mass ratio in the restricted three-body problem.
Since all orbits are symmetric with respect to the anti-symplectic involution
\begin{align}\label{eq reflection}
	\rho \colon (q_1, q_2, q_3, p_1, p_2, p_3) \mapsto (q_1, -q_2, q_3, -p_1, p_2, p_3),
\end{align}
they can be found by a perpendicular shooting method.

%% file: hypersurfaces.tex
The main goal of this chapter is to compute generators of the first de Rham cohomology, which is essential to us for the obstruction to contact forms by closed orbits.
These generators are found by computing the fundamental group of the Moser-regularised energy hypersurface, then abelianising it to the first homology group, modding out torsion to get real coefficients and, finally, dualising to get generators of the first de Rham cohomology.
We will write all groups that appear here multiplicatively unless they are inherently abelian, in which case we will write them additively.
First of all, we will compute the planar case which will take most of this chapter and then comment on the spatial case.

Recall the notation~$\Sigma_c := H^{-1}(c)$ for the energy hypersurface of the Kepler problem or the restricted three-body problem.
Denote by~$\overline{\Sigma}_c$ the corresponding regularised energy hypersurface, where collisions have been added by Moser regularisation.
In the first step we compute the fundamental group of~$\overline{\Sigma}_c$ using the well-known theorem of Seifert-van Kampen.

We will explicitly compute the fundamental group of the bounded component of the regularised energy hypersurface of the Kepler problem and then use the relations found there to compute the more complicated hypersurface of the restricted three-body problem above the highest critical value~$H_\mu(L_5)$.

In the Moser regularisation first the roles of~$P$ and~$Q$ are interchanged, such that the base points of the cotangent bundle now corresponded to the momentum of the particle and the fibre to its position.
At every point in the base the intersection between the hypersurface and the fibre is then a circle of positions.
The base as points of momentum is endowed with the metric of the stereographic projection of $S^2$ through the north pole~$\mathcal{N}$, corresponding to infinite momentum at collision.
The regularised energy hypersurface is thus the unit cotangent bundle~$S^*(S^2) \cong \RR \PP^3$ of the round 2-sphere.
We choose as the second chart of~$S^2$ the stereographic projection through the south pole~$\mathcal{S}$, corresponding to zero momentum.
So, we define the subsets
\begin{align*}
	U_1 &:= S^*(S^2 \setminus \{ \mathcal{N} \}) \cong S(\RR^2)
	\\
	U_2 &:= S^*(S^2 \setminus \{ \mathcal{S} \}) \cong S(\RR^2) \text{ and}
	\\
	U_3 &= S^*(S^2 \setminus\{ \mathcal{S}, \mathcal{N} \}) \cong S(\RR^2 \setminus \{ 0 \}),
\end{align*}
where the trivialisations of~$U_1$ and~$U_3$ are given by the stereographic projection through the north pole and the trivialisation of~$U_2$ by the stereographic projection through the south pole.
The change of these variables is given in local coordinates~$x \in \RR^2$ by
\begin{align*}
	\Phi(x_1,x_2) = \left( \frac{x_1}{x_1^2+x_2^2}, \frac{x_2}{x_1^2+x_2^2} \right)
\end{align*}
and the Jaconian is
\begin{align*}
	\Diff \Phi(x_1,x_2) = \frac{1}{(x_1^2+x_2^2)^2}
	\begin{pmatrix}
		x_2^2-x_1^2		&	-2x_1 x_2
		\\
		-2 x_1 x_2		&	x_1^2 - x_2^2
	\end{pmatrix}.
\end{align*}
As the base point for the fundamental groups we choose
\begin{align*}
	x_0 := ((1,0),(1,0)) \in S(\RR^2 \setminus \{ 0 \}) \cong U_3 = U_1 \cap U_2,
\end{align*}
\ie the point with momentum $P=(1,0)$ and position~$Q=(1,0)$.
In the trivialisation of~$U_2$ this point~$x_0$ corresponds to $((1,0),(-1,0)) \in S(\RR^2)$.

The fundamental groups of the subsets are
\begin{align*}
	\pi_1(U_1, x_0) &= \langle \zeta \rangle \cong \ZZ
	\\
	\pi_1(U_2, x_0) &= \langle \eta \rangle \cong \ZZ
	\\
	\pi_1(U_3, x_0) &= \langle \xi_1, \xi_2 \mid [\xi_1, \xi_2] \rangle \cong \ZZ \times \ZZ,
\end{align*}
where~$\zeta$ and $\xi_1$ is each the class of homotopic loops based at~$x_0$ and represented by
$$
t \mapsto ((1,0),(\cos(2\pi t), \sin(2 \pi t))),
$$
\ie a simple loop in the position fibre,
$\eta$ is also represented by a loop
$$
t \mapsto ((1,0),(-\cos(2\pi t), \sin(2 \pi t))),
$$
in the fibre and $\xi_2$ is represented by a loop
$$
t \mapsto ((\cos (2 \pi t), \sin (2 \pi t)),(1,0))
$$
in the momentum base.

For the theorem of Seifert-van Kampen we need to compute the images~$v_i(\xi_j)$ for $i,j \in \{1,2\}$ of generators of~$\pi_1(U_3,x_0)$ after the homomorphisms~$v_i$ induced by the inclusions $U_i \hookrightarrow U_3$.
Since we chose the same trivialisation for~$U_1$ and~$U_3$, the first two are simply~$v_1(\xi_1) = \zeta$ and~$v_1(\xi_2)=1$ because~$\zeta$ and~$\xi_1$ are identically represented in the trivialisation and the representation of $\xi_2$ is contractible in~$U_1$.
For the second homomorphism~$v_2$ we need to check the differential of the change of trivialisations, \ie the Jacobian of the change of coordinates between the stereographic projection through the north and the south pole.

The representation of~$\xi_1$ readily gets mapped onto the representation of~$\eta$, so~$v_2(\xi_1) = \eta$.
The image of the representation of~$\xi_2$ is again contractible in the base, but it twists the fibre twice in the opposite direction of~$\eta$ since
\begin{align*}
	\Diff \Phi(\cos( 2 \pi t), \sin (2 \pi t))
	\begin{pmatrix}
		1 \\ 0
	\end{pmatrix}
	&=
	\left( \sin^2( 2 \pi t) - \cos^2( 2 \pi t), - 2 \sin( 2 \pi t) \cos( 2 \pi t) \right)
	\\
	&=
	\left( - \cos( 4 \pi t), - \sin( 4 \pi t) \right).
\end{align*}
All in all we have
\begin{align*}
	v_1(\xi_1) &= \zeta,
	&
	v_1(\xi_2) &= 1,
	\\
	v_2(\xi_1) &= \eta,
	&
	v_2(\xi_2) &= \eta^{-2}
\end{align*}
and the fundamental group of~$\overline{\Sigma}_c$ becomes
\begin{align*}
	\pi_1(\overline{\Sigma}_c,x_0)
	=
	\langle \zeta, \eta \mid \zeta = \eta, \eta^{-2} \rangle
	=
	\langle \zeta \mid \zeta^2 \rangle \cong \ZZ_2.
\end{align*}
Of course we would have known that earlier from the fact that~$\overline{\Sigma}_c \cong S^*S^2 \cong \RR \PP^3$, but now we can use the relations from above to compute the fundamental groups of more complicated regularised hypersurfaces.

The surface we are interested in is the energy level set of the restricted three-body problem above the highest critical value~$H_\mu(L_5) < c$.
Here, we need to regularise two singularities and we will therefore apply the theorem of Seifert-van Kampen twice.
As the subsets we again choose the unregularised energy level set
$$
U_1 := \Sigma_c \cong S^*(\RR^2 \setminus \{M_1, M_2\}) \cong S (\RR^2 \setminus \{M_1, M_2\}),
$$
for the local regularising charts each a copy of the unit cotangent bundle of a small open 2-disc
$$
U_2,U_3 := S^*(B_\varepsilon) \cong S(\RR^2),
$$
and the intersections become unit cotangent bundles of punctured 2-discs
$$
U_4, U_5 := S^*(B_\varepsilon \setminus \{ 0 \} ) \cong S(\RR^2 \setminus \{0\}).
$$
Remember, however, that the trivialisation of~$\Sigma_c$, where currently the position coordinates form the base, is changed in the first step of regularisation, such that the momentum becomes the base and the fibre is the position.
The fundamental groups of these spaces are
\begin{align*}
	\pi_1(U_1) &= \langle r, w_1, w_2 \mid [r,w_1], [r,w_2] \rangle,
	\\
	\pi_1(U_2) &= \langle \eta_2 \rangle,
	\\
	\pi_1(U_3) &= \langle \eta_3 \rangle,
	\\
	\pi_1(U_4) &= \langle \xi_1, \xi_2 \mid [\xi_1, \xi_2] \rangle
	\quad \text{and}
	\\
	\pi_1(U_5) &= \langle \xi_3, \xi_4 \mid [\xi_3, \xi_4] \rangle,
\end{align*}
where~$r$ is represented by the loop in momentum coordinates over a fixed position, $w_1$ is the winding in position around~$M_1$ and~$w_2$ is the winding in position around~$M_2$, both commuting with~$r$.
The two~$\eta_2$ and~$\eta_3$ are defined, same as~$\eta$ above, as the loop in position coordinates, as well as $\xi_1$ and~$\xi_3$, just as in~$\xi_1$ from above, while~$\xi_2$ and~$\xi_4$ are a represented by a loop in momentum coordinates, as was~$\xi_2$ from before.
Denote the homomorphisms of fundamental groups induced by inclusions as
\begin{align*}
	\begin{split}
		v_1 &\colon \pi_1(U_4) \to \pi_1(U_1),
		\\
		v_2 &\colon \pi_1(U_4) \to \pi_1(U_2),
	\end{split}
	\begin{split}
		v_3 &\colon \pi_1(U_5) \to \pi_1(U_1),
		\\
		v_4 &\colon \pi_1(U_5) \to \pi_1(U_3).
	\end{split}
\end{align*}
We can now use the same relations as in the Kepler problem:
\begin{align*}
	v_1(\xi_1) &= w_1
	&
	v_1(\xi_2) &= r
	\\
	v_2(\xi_1) &= \eta_2
	&
	v_2(\xi_2) &= \eta_2^{-2}
	\\
	v_3(\xi_3) &= w_2
	&
	v_3(\xi_4) &= r
	\\
	v_4(\xi_3) &= \eta_3
	&
	v_4(\xi_4) &= \eta_3^{-2}
\end{align*}
The only difference is that the loop in momentum coordinates is no longer contractible in the original~$\Sigma_c = U_1$.
Putting together the generators and relations, we get
\begin{align*}
	\pi_1(\overline{\Sigma}_c)
	&= \langle r, w_1, w_2, \eta_2, \eta_3 \mid [r,w_1], [r,w_2], w_1 = \eta_2, r = \eta_2^{-2}, w_2 = \eta_3, r = \eta_3^{-2} \rangle
	\\
	&= \langle r, w_1, w_2 \mid [r,w_1], [r,w_2], r = w_1^{-2}, r= w_2^{-2} \rangle
	\numberthis \label{eq relations fundamental}
	\\
	&= \langle w_1, w_2 \mid w_1^2 = w_2^2 \rangle.
\end{align*}

In order to find the first de Rham cohomology, we first compute the abelianisation 
\begin{align*}
	H_1(\overline{\Sigma}_c, \ZZ) \cong \pi_1(\overline{\Sigma}_c)^\mathrm{ab}
	&= \langle w_1, w_2 \mid 2 w_1 = 2 w_2 \rangle
	\\
	&= \langle w_1, w_2, z \mid 2 w_1 = 2 w_2, z = w_2 - w_1 \rangle
	\\
	&= \langle w_1, z \mid 2 z \rangle
	\\
	&\cong \ZZ \times \ZZ_2,
\end{align*}
which is isomorphic to the first homology group with integer coefficients.
By modding out torsion by the universal coefficient theorem we get the first homology with real coefficients
\begin{align*}
	H_1(\overline{\Sigma}_c, \RR) \cong \RR \cong H^1_\mathrm{dR}(\overline{\Sigma}_c),
\end{align*}
which is also isomorphic to the first de Rham cohomology by duality and de Rham's theorem.
We denote the push forward of the inclusion~$\iota \colon \Sigma_c \hookrightarrow \overline{\Sigma}_c$ on homology by
\begin{align*}
	\iota_* \colon H_1(\Sigma_c) \cong \RR^3 &\to H_1(\overline{\Sigma}_c) \cong \RR.
\end{align*}
The kernel of~$\iota_*$ can be found from line~\eqref{eq relations fundamental} to be~$\ker(\iota_*) = \langle r + 2w_1, r+2w_2 \rangle_\RR$ and the coimage is then~$\ker(\iota_*)^\perp = \langle  2 r - w_1  -w_2 \rangle_\RR$, where we abuse the notation from above of fundamental classes~$r$, $w_1$ and~$w_2$ to also denote generators of homology.

Dualising to cohomology, we see that the image of the pullback of~$\iota$ on cohomology
\begin{align*}
	\iota^* \colon H_\mathrm{dR}^1 (\overline{\Sigma}_c) \cong \RR \to H_\mathrm{dR}^1(\Sigma_c) \cong \RR^3
\end{align*}
is then $\mathrm{im}(\iota^*) = \langle 2 \diff \vartheta - \diff \varphi_1 - \diff \varphi_2 \rangle$, where~$\vartheta$ is the polar angle in momentum coordinates and~$\varphi_1$ and~$\varphi_2$ are polar angles in position coordinates centred at~$M_1$ and~$M_2$, respectively.
So, we have found a generator of the first de Rham cohomology of the regularised energy hypersurface above the highest critical value, which we can compute easily for periodic non-collision orbits by twice the rotation number minus the two winding numbers around the primaries.
We summarise the results from this chapter on the planar restricted three-body problem in the following lemma:

\begin{lemma} \label{lemma de Rham generator}
	For the planar circular restricted three-body problem and energies~$c$ above the highest critical value~$H_\mu(L_5)$ the first de Rham cohomology of the Moser-regularised energy hypersurface~$\overline{\Sigma}_c$ is one-dimensional and has generator~$0 \neq [ \beta_0] \in H_\mathrm{dR}^1(\overline{\Sigma}_c)$ which agrees with~$2 \diff \vartheta - \diff \varphi_1 - \diff \varphi_2$ on the unregularised level set~$\Sigma_c$.
\end{lemma}

If we use the same setup in the spatial case, we see that all sets
\begin{align*}
	U_1
	&=
	\Sigma_c \cong S(\RR^2 \setminus \{M_1, M_2\}),
	\\
	U_2,U_3
	&=
	S^*(B_\varepsilon) \cong S(\RR^2) \text{ and}
	\\
	U_4, U_5
	&=
	S^*(B_\varepsilon \setminus \{ 0 \} ) \cong S(\RR^2 \setminus \{0\})
\end{align*}
are simply connected.
Consequently, by the theorem of Seifert-van Kampen also the union~$\overline{\Sigma}_c$ is simply connected.
This means that in the spatial restricted three-body problem every closed orbit is contractible in the regularised as well as in the unregularised energy level set.

%% file: generating_orbits.tex
Next, we will describe the notion of generating orbits for the restricted three-body problem.
We will introduce notation but also describe some results from other authors and connect their results to our purposes.

Generating orbits are limits of orbits of the restricted three-body problem.
In our setting we will use the limit as~$\mu \to 0$ but in general also other limits, for example letting the angular momentum tend to zero, might be possible.
This chapter establishes terminology and results around these generating orbits which is taken from~\cite{Henon1}.
Another extensive work on generating orbits is~\cite{Bruno}, where more theory is explained and mainly rotating coordinates are used.
In~\cite{Henon1} fixed coordinates are used to describe the generating orbits which makes it easier for us to use the geometry of Kepler ellipses and compute the action of generating orbits.
We are only interested in periodic orbits here and, hence, we will only consider periodic generating orbits.

\begin{definition}\label{def generating orbit}
	Let~$\gamma_\mu$ be a periodic orbit of the restricted three-body problem with mass ratio~$\mu>0$.
	Then~$\gamma$ is called a \emph{generating orbit} if there exists a sequence of orbits~$\gamma_\mu$ such that $\gamma_\mu \to \gamma$ as~$\mu \to 0$.
\end{definition}

\begin{remark}
	In general, generating orbits are not orbits of the restricted three-body problem for~$\mu=0$, \ie the rotating Kepler problem, and vice-versa.
	Furthermore, when we use the notion of a generating orbit, we usually mean this in Henon's terms, who mostly worked with numerical methods.
	For an analytical proof of the fact that there exist continued orbits in the restricted three-body problem we rely on section~\ref{section continuation of generating orbits}.
\end{remark}

Next, we will classify \emph{species} of generating orbits.
This notion goes back to Poincaré in~\cite{Poincare3} who called his predicted periodic orbits with near collisions in the general three-body problem ``solutions périodiques de deuxième espèce''.
We will adopt Hénon's more general definition of orbit species:
\begin{definition}\label{def species}
	A generating orbit is of the \textit{first species} if it is a Keplerian orbit, it is of the \textit{second species} if at least one point coincides with~$M_2$ and it is of the \textit{third species} if it only consists of~$M_2$.
\end{definition}
This definition does not give mutually exclusive species.
In fact, third species generating orbits are always also of both the first and the second species.
There are second species orbits that are also of the first species but there are also orbits that are exclusively of the first or second species.
However, this definition is in terms of Hénon's \emph{principle of positive definition} where ``a definition relating to orbits in a family should not be based on a negative property, such as an inequality''.
This principle gives families of generating orbits the same species at the cost of exclusivity of species.

\subsection{First species} \label{section first species}
First species orbits are periodic orbits of the rotating Kepler problem.
We define in accordance with the principle of positive definition:

\begin{definition}
	A generating orbit of the first species is called of the \emph{first kind} if it is a circular orbit and of the \emph{second kind} if in fixed coordinates $M_2$ and~$M_3$ each make an integral number of revolutions~$I$ and~$J$ around~$M_1$.
\end{definition} 
Again, this definition is not exclusive and orbits belonging to both kinds are bifurcation orbits.
More details on the first species can be found in~\cite{Henon1}.

\subsection{Second species} \label{section second species}
According to definition~\ref{def species} a second species orbit passes through~$M_2$ at least once.
We call this event a~\emph{collision}.
A periodic second species generating orbit will collide infinitely many times, but there can also be multiple collisions during one minimal period.
Abiding by the principle of positive definition, we declare a finite piece of a Keplerian orbit which begins and ends in collision an \emph{arc}.
Note that also an arc can include collisions, subdividing the arc into \emph{basic arcs}.
The angle at collision between basic arcs is called the \emph{deflection angle} and a generating orbit of the second species is called \emph{ordinary generating orbit} if all deflection angles are nonzero.
Non-ordinary generating orbits are again bifurcation orbits.

The Kepler orbit belonging to the arc is called the \emph{supporting Kepler orbit} or the \emph{supporting Kepler ellipse}, if referring to the geometrical object.
Each second species generating orbit consists of a sequence of arcs $U_1, U_2, \dots U_k$, which is repeated periodically, and each arc is fully defined by its supporting Kepler solution and times~$t'_j$ and~$t''_j$ of initial and final collision.
The \emph{duration} of an arc~$U_j$ is~$\tau_j := t''_j - t'_j$.

The further study of second species orbits will from now on be confined to the study of arcs where we will continue to only state results and definitions which we will work with later.

\begin{figure}[tb]
	\centering
	\begin{subfigure}[t]{0.32\textwidth}
		\centering
		\includegraphics[width=\textwidth]{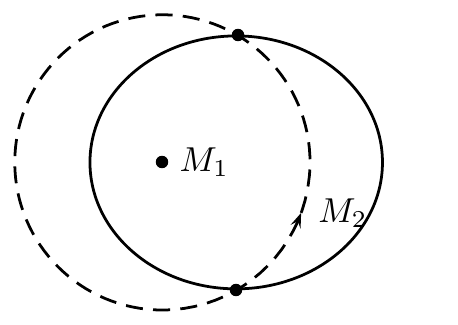}
		\caption{Type 1: transverse intersection.}
	\end{subfigure}
	\hfill
	\begin{subfigure}[t]{0.32\textwidth}
		\centering
		\includegraphics[width=\textwidth]{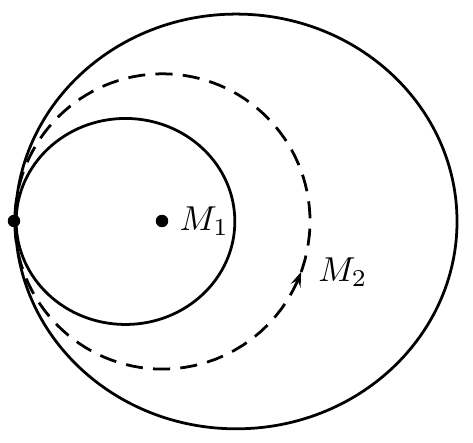}
		\caption{Type 2: tangent to unit circle.}
	\end{subfigure}
	\hfill
	\begin{subfigure}[t]{0.32\textwidth}
		\centering
		\includegraphics[width=\textwidth]{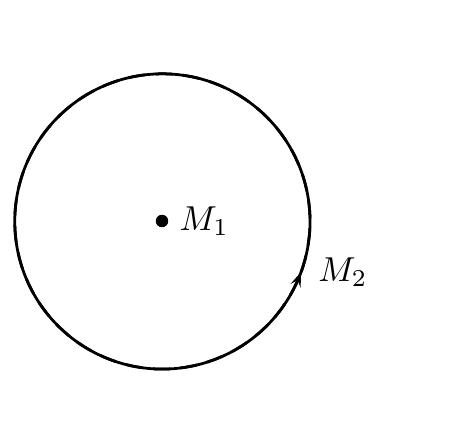}
		\caption{Type 3 and 4: identical with unit circle in retrograde and direct direction.}
	\end{subfigure}
	\caption{Types of second species supporting ellipses.}
\end{figure}
Let~$r_1$ be the pericentre distance and~$r_2$ the apocentre distance of the supporting ellipse.
For a collision to occur, we need~$r_1 \leq 1 \leq r_2$.
We will distinguish the following cases:
For $r_1<1<r_2$ we have a non-circular supporting ellipse which intersects the unit circle in two distinct points transversally.
The corresponding arcs will be called of \emph{type 1}.
For $r_1=1$ and~$r_2>1$ or~$r_1<1$ and~$r_2=1$ the supporting Kepler ellipse is tangent to the unit circle and we will call the arcs of \emph{type 2}.
In the remaining case $r_1=r_2=1$ the supporting Kepler ellipse is identical with the unit circle and we will call the corresponding arc of \emph{type 3} if it is retrograde and of \emph{type 4} if it is direct.
Type 1 is the most interesting and comes in families while the other types are isolated.

We will further subdivide type 1 arcs into \emph{$S$-arcs} which in fixed coordinates begin and end at different points on the unit circle and \emph{$T$-arcs} which begin and end at the same point.
A type~1 arc will furthermore be called \emph{ingoing} if at the initial collision its velocity vector points to the inside of the unit circle, and \emph{outgoing} else.
Both $S$-arcs and $T$-arcs can be ingoing and outgoing.
\subsubsection{S-arcs}
An $S$-arcs is symmetric with respect to the major axis of the supporting ellipse in fixed coordinates and intersects this axis $2 J + 1$ times for~$J \geq 0$.
Let~$R$ be the central---\ie the~$J+1$\textsuperscript{st}---intersection point.
Then $R$ lies either at the pericentre or apocentre and we call it the \emph{midpoint}.
Since both collision points lie at~$(1,0)$ in rotating coordinates, the midpoint~$R$ also lies on the $q_1$-axis and the arc is symmetric with respect to~$\rho$ from~\eqref{eq reflection}.

\subsubsection{T-arcs}
T-arcs begin and end in the same point in fixed coordinates and are therefore full Kepler ellipses.
Analogously to rotating Kepler orbits, the semi-major axis of the supporting ellipse can be expressed by numbers~$I$ and~$J$ of rotation of~$M_2$ and~$M_3$ around~$M_1$ by $a = \left( I /J \right)^{2/3}$.
Because we have
\begin{proposition}[proposition 4.3.2 in \cite{Henon1}]\label{prop no consecutive T-arcs}
	An ordinary generating orbit of the second species can not contain two identical $T$-arcs of type~1 in succession.
\end{proposition}
the numbers~$I$ and~$J$ must again be relatively prime, \ie no multiple covers are allowed.
$T$-arcs are not symmetric and for every mutually prime~$I$ and~$J$ there exists a family~$T_{I,J}^i$ of ingoing $T$-arcs and a family~$T_{I,J}^e$ of outgoing $T$-arcs.
\begin{figure}[tb]
	\centering
	\begin{subfigure}[t]{0.495\textwidth}
		\centering
		\includegraphics[width=\textwidth]{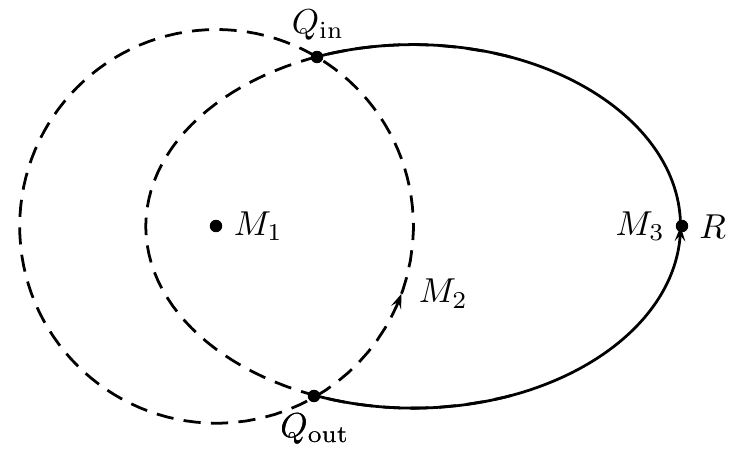}
	\end{subfigure}
	\hfill
	\begin{subfigure}[t]{0.495\textwidth}
		\centering
		\includegraphics[width=\textwidth]{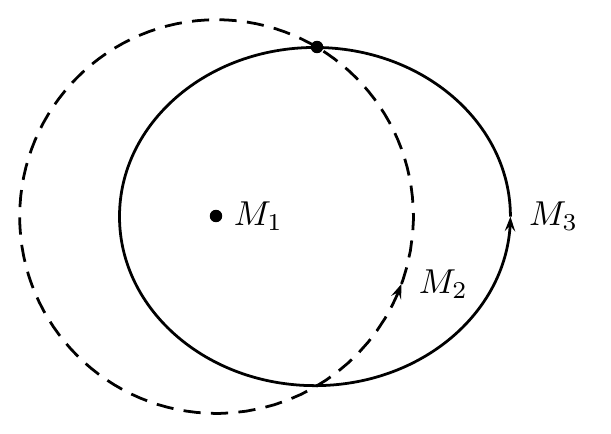}
	\end{subfigure}
	\caption{An outgoing S-arc and an ingoing T-arc.}
\end{figure}

\subsection{Continuation of second species generating orbits} \label{section continuation of generating orbits}

In order to show analytically that the orbits described above are actually generating orbits we need to find converging sequences of orbits in the restricted three-body problem as~$\mu \to 0$.
Conversely, we then have for every generating orbit an arbitrarily close orbit in the restricted problem for some small enough mass ratio.

For first species orbits there are many works on the existence of such sequences, for example~\cite{Poincare1}, \cite{Birkhoff} and \cite{hagihara} for the first kind, \cite{Arenstorf} for symmetric, and~\cite{Bruno} for asymmetric second kind orbits.
The remainder of his chapter is focused on second species generating orbits.

\subsubsection{A more general result}\label{subsection sketch of proof of BMK}
The main theorem from~\cite{BolotinMacKay} helps to show that many ordinary second species orbits are actually generating orbits.
The setting is somewhat more general in that paper and we can also follow from the proof presented there that the action of the orbits in the restricted three-body problem converges to the action of the generating orbit.
We will therefore first describe the notation which we adapt for our purposes and then explain in more detail the result and the relevance for the present work.

Let~$\mathcal{Q}$ be a two- or three-dimensional smooth manifold and~$\mathcal{P}=\{p_1, \dots p_n\} \subset \mathcal{Q}$ a finite subset.
The set~$\mathcal{Q} \setminus \mathcal{P}$ shall be the configuration space for the Hamiltonian
\begin{align*}
	H_\mu(q,p) = H_0(q,p) + \mu V(q),
\end{align*}
where $H_0(q,p) = \frac{1}{2} \norm{p + A_q}^2_{g^*_q} + W(q)$ is a magnetic Hamiltonian and~$V$ is another smooth potential with Newtonian singularities at every~$p_i \in \mathcal{P}$.
So, the Hamiltonian~$H_0$ is defined on all of~$T^*\mathcal{Q}$, while~$H_\mu$ is defined on~$T^*(\mathcal{Q} \setminus \mathcal{P})$.
For the planar restricted three-body problem with small mass ratio~$\mu$ we choose~$\mathcal{Q} = \RR^2 \setminus \{0\}$ and~$\mathcal{P}=\{(1,0)\}$, \ie we shift the coordinates such that the origin is always at the heavier primary~$M_1$.
The Hamiltonian then is of the form from above, with
\begin{align*}
	A_q &= q_2 \diff q_1 - q_1 \diff q_2,
	\\
	W(q) &= - \frac{1}{\norm{q}} - \frac{\norm{q^2}}{2}
	\qquad \text{and}
	\qquad
	V(q) = \frac{1}{\norm{q}} - \frac{1}{\norm{q - (1,0)}} - q_1.
\end{align*}

Fix an energy~$c$ such that the open Hill's region~$\mathfrak{H}_c := \{q \in \mathcal{Q} \mid W(q) < c \}$ contains~$\mathcal{P}$.
Suppose we have a finite set~$K$ of nondegenerate collision orbits~$\gamma_k \colon [0, \tau_k] \to \mathfrak{H}_c$ of the Hamiltonian system~$H_0$ such that $\gamma_k(0)=p_{\alpha_k}$, $\gamma_k(\tau_k) = p_{\beta_k} \in \mathcal{P}$ and~$\gamma(t) \in \mathfrak{H}_c \setminus \mathcal{P}$ for all other~$t \in (0,\tau_k)$.
A~\emph{chain} is a sequence~$(\gamma_{k_i})_{i \in \ZZ}$ of orbits in~$K$ such that additionally~$\gamma_{k_i}(\tau_{k_i}) = \gamma_{k_{i+1}}(0)$ and~$\dot{\gamma}_{k_i}(\tau_{k_i}) \neq \dot{\gamma}_{k_{i+1}}(0)$, \ie they are connected collision orbits with nonzero deflection angle.
Let~$\mathcal{W}_k$ be open neighbourhoods for each of the sets~$\gamma_k([0,\tau_k])$ in~$\mathcal{Q}$.
An orbit~$\gamma \colon \RR \to \mathfrak{H}_c$ is said to \emph{shadow} the chain~$(\gamma_{k_i})_{i \in \ZZ}$ if there exists an increasing sequence~$(t_i)_{i \in \ZZ}$ such that $\gamma([t_i, t_{i+1}]) \subset \mathcal{W}_{k_i}$.

The \emph{nondegeneracy} of such orbits~$\gamma$ is defined as the Morse-nondegeneracy of the critical point $(u,\tau) \in W^{1,2}(p_\alpha, p_\beta) \times \RR^+$ of the action functional~$\mathcal{A}(\gamma)$, where~$W^{1,2}(p_\alpha, p_\beta)$ is the space of all $W^{1,2}$-functions $u \colon [0,1] \to \mathcal{Q}$ with fixed endpoints~$u(0) = p_\alpha$, $u(1) = p_\beta$ and $\gamma (t) = u(t/\tau)$.
There are four other equivalent ways to define the nondegeneracy described in~\cite{BolotinMacKay}.
The only other one we will use here is the following:
Denote by~$q(\lambda, t)$ the general solution of Hamilton's second order differential equations of motion in the configuration space with parameter~$\lambda$ and by $h(\lambda)$ the Hamiltonian energy.
Then the collision orbit is nondegenerate if the system
\begin{align}\label{eq nondegeneracy of orbits}
	q(\lambda, 0) = p_\alpha,
	\quad
	q(\lambda, \tau) = p_\beta,
	\quad
	h(\lambda) = c
\end{align}
has full rank $2 \dim \mathcal{Q} + 1$.

\begin{theorem}[theorem 1.1 from~\cite{BolotinMacKay}] \label{theorem Bolotin MacKay}
	There exists~$\mu_0 > 0$ such that for all $\mu \in (0, \mu_0]$ and any chain $(\gamma_{k_i})_{i \in \ZZ}$ the following holds:
	\begin{itemize}
		\item
		There exists a trajectory $\gamma \colon \RR \to \mathfrak{H}_c$ of energy $c$ for the system $H_\mu$ shadowing the chain $(\gamma_{k_i})_{i \in \ZZ}$ and it is unique up to a time-shift if the neighbourhoods~$\mathcal{W}_k$ are chosen small enough.
		
		\item
		The orbit $\gamma$ converges to the chain of collision orbits as $\mu \to 0$:
		There exists an increasing sequence $(t_i)_{i \in \ZZ}$ such that
		\begin{align*}
			\max_{t_i \leq t \leq t_{i+1}} \mathrm{dist}(\gamma(t), \gamma_{k_i}([0,\tau_{k_i}])) \leq \mu C_1,
		\end{align*}
		where the constant~$C_1>0$ depends only on the set~$K$ of collision orbits.
		
		\item
		If we additionally have $\dot{\gamma}_{k_i}(\tau_{k_i}) \neq - \dot{\gamma}_{k_{i+1}}(0)$, the orbit~$\gamma$ avoids collision by a distance of order~$\mu$:
		there exists a constant~$C_2 \in (0, C_1)$, depending only on~$K$, such that
		\begin{align*}
			\mu C_2 \leq \min_{t_{i-1} \leq t \leq t_{i+1}} \mathrm{dist} (\gamma(t) , p_{\alpha_{k_i}}).
		\end{align*}
	\end{itemize}
\end{theorem}

\subsubsection{Application to the restricted three-body problem}\label{subsection application of BMK to PCR3BP}

We now want to apply this theorem~\ref{theorem Bolotin MacKay} to the restricted three-body problem and find orbits for positive mass ratios~$\mu>0$ shadowing chains of collision orbits, \ie second species generating orbits.
In the same paper~\cite{BolotinMacKay} a partial result is already shown:

\begin{lemma}[lemma 3.2 in \cite{BolotinMacKay}] \label{lemma non-degenerate T-arcs}
	For $c \in (-3/2, \sqrt{2})$ the collision orbit at~$\mu = 0$ in the rotating frame corresponding to a whole number~$I$ of revolutions of an ellipse with rational frequency $I/J \in A_c$ in the allowed set of frequencies~$A_c$ for energy~$c$ starting and ending at a collision with~$M_2$ is nondegenerate.
\end{lemma}

In the language of generating orbits this means that every ordinary generating orbit composed of T-arcs is indeed a generating orbit.
We would like to use generating orbits composed of S-arcs and we will therefore have to prove accordingly:

\begin{lemma}\label{lemma non-degenerate S-arcs}
	For $c \in (-3/2, \sqrt{2})$ a non-rectilinear ordinary S-arc collision orbit starting and ending in~$M_2$ with semi-major axis~$a>1$ of the supporting ellipse is nondegenerate.
\end{lemma}

\begin{proof}
	This proof is very similar to the proof of lemma~\ref{lemma non-degenerate T-arcs} in~\cite{BolotinMacKay} but a lot harder to actually compute the estimates.
	Enforcing the first equation of~\eqref{eq nondegeneracy of orbits}, we have $\dim \mathcal{Q} = 2$ free parameters left for~$\lambda$.
	We parametrize the supporting ellipse through our initial point of collision~$Q_0$ by the position of the second focus~$F \in \RR^2$ using two variables:
	the polar angle~$\theta$ between~$Q_0$ and~$F$ and the semi-major axis
	\begin{align*}
		a
		&= \frac{1}{2}(1 + \mathrm{dist}(Q_0,F))
		\\
		&= \frac{1}{2}(1 + \sqrt{1 + \mathrm{dist}(M_1, F)^2 - 2 \cos (\theta) \mathrm{dist}(M_1,F)}).
	\end{align*}
	This is a good parameter if $\mathrm{dist}(M_1,F) \neq \cos \theta$ or, equivalently, if $\mathrm{dist}(Q_0,F) \neq \sin \theta$, so we restrict to $\mathrm{dist}(Q_0,F)> 1$, \ie if and only if $a>1$.
	We will see later on that this does effectively not restrict S-arcs in the direct sense of rotation~$\epsilon' = +1$, which are the only ones we need in the main proof.
	\begin{figure}[tb]
		\centering
		\includegraphics{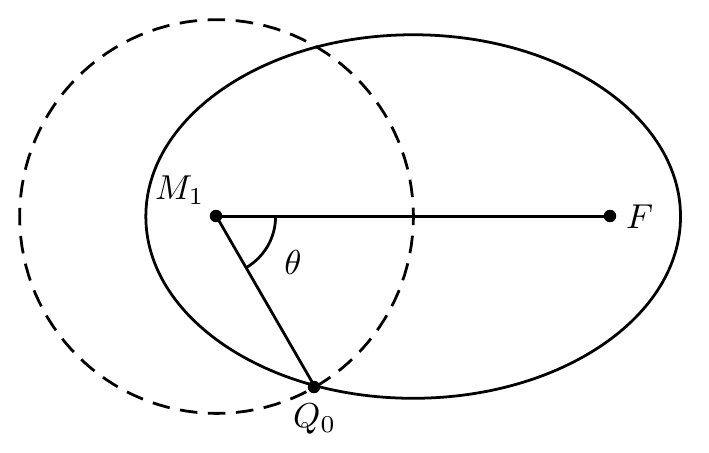}
		\caption{Parametrisation of second species ellipses.}
	\end{figure}
	
	The second condition of~\eqref{eq nondegeneracy of orbits} is satisfied nondegenerately since the supporting ellipse intersects the unit circle transversely and both~$M_2$ and~$M_3$ have nonzero velocities.
	Fixing the endpoint also fixes the angle~$\theta$, so the remaining free parameter is~$a$.
	
	We are left to show that the derivative of the energy by~$a$ is nonzero.
	Dependent on~$a$ and~$\theta$ we compute the eccentricity and the energy, as will be also done in section~\ref{subsection fixed theta} later on:
	\begin{align*}
		\epsilon
		&=
		\frac{\cos \theta + \sqrt{4a^2 - 4a + \cos^2 \theta}}{2a}
		\\
		H_0
		&=
		- \frac{1}{2a} - \epsilon' \sqrt{ 1 - \cos \theta  \frac{\cos \theta + \sqrt{4a^2 - 4a  + \cos^2 \theta}}{2a}}
	\end{align*}
	We exclude the case of rectilinear orbits with $\theta=0$ due to the root being zero and the sign of~$\epsilon'$ changing here.
	For numbers of revolution~$J>0$ there would also occur a collision with~$M_3$ and~$M_1$ in that case which one would have to deal with additionally.
	
	For all other orbits with~$\theta \in (0,\pi)$ and~$a>1$ we differentiate~$H_0$ by~$a$ to get
	\begin{align*}
		\frac{\partial H_0}{\partial a}
		&=
		\frac{1}{2a^2} + \epsilon' \frac{\cos \theta \frac{\frac{4a^2-2a}{\sqrt{4a^2-4a+\cos^2 \theta}}-\cos \theta - \sqrt{4a^2 - 4a + \cos^2 \theta}}{2a^2}}{2 \sqrt{1- \cos \theta \frac{\cos \theta + \sqrt{4a^2 - 4a + \cos^2 \theta}}{2a}}}.
	\end{align*}
	This vanishes if and only if
	\begin{align*}
		2 \sqrt{1- \cos \theta \frac{\cos \theta + \sqrt{4a^2 - 4a + \cos^2 \theta}}{2a}}
		=
		\hspace{-9em}
		\\
		\qquad &=
		-\epsilon'\cos \theta \left( \frac{4a^2-2a}{\sqrt{4a^2-4a+\cos^2 \theta}}-\cos \theta - \sqrt{4a^2 - 4a + \cos^2 \theta} \right)
		\\
		&=
		\frac{-2\epsilon' a \cos \theta}{\sqrt{4a^2-4a+\cos^2 \theta}}\underbrace{\left( 1- \cos \theta \frac{\cos \theta + \sqrt{4a^2-4a+\cos^2 \theta}}{2a}\right)}_{=a(1-\epsilon^2)>0} .
	\end{align*}
	Based on this equation, we can distinguish four cases:
	Case 1 is $\epsilon'=+1$ and $\cos \theta \geq 0$, case 2 is $\epsilon'=+1$ and $\cos \theta < 0$, case 3 is $\epsilon'=-1$ and $\cos \theta > 0$, and case 4 is $\epsilon'=-1$ and $\cos \theta \leq 0$.
	
	In cases 1 and 4 the right-hand side becomes nonpositive and the nondegeneracy is obvious since the left-hand side is always positive for~$\epsilon \neq 1$.
	For the remaining cases we will have to work a little bit harder.
	In the following computations we will denote the left-hand side by~$u=u(a,\theta)$ and the right-hand side by~$v=v(a,\theta)$.
	
	In case 2 the sign of~$\cos \theta$ is negative, so we can estimate the left-hand side by $u>2$.
	To show that the right-hand side~$v$ is smaller, we insert the boundary value~$a=1$ to get
	\begin{align*}
		v(1, \theta)
		&=
		\frac{-2 \cos \theta}{\betrag{\cos \theta}}\left( 1- \cos \theta \left( \cos \theta + \betrag{\cos \theta} \right) \right)
		=
		2
	\end{align*}
	and then see that the derivate~$\partial v / \partial a$ is negative by
	\begin{align*}
		\frac{\partial v}{\partial a}
		&=
		\\
		&
		\hspace{-20pt}
		- \cos \theta
		\left(
		\frac{
			(8a-2)\sqrt{4a^2-4a+\cos^2 \theta}-(4a^2-2a)\frac{8a-4}{2\sqrt{4a^2-4a+\cos^2\theta}}
		}{4a^2-4a + \cos^2 \theta}
		-
		\frac{8a-4}{2\sqrt{4a^2-4a+\cos^2\theta}}
		\right)
		\\
		&=
		- \cos \theta \frac{4a \cos^2 \theta - 4a}{(4a^2-4a+\cos^2 \theta)^{\frac{3}{2}}}<0.
	\end{align*}
	
	In the last remaining case 3 we have~$\cos \theta > 0$ and~$\epsilon'=-1$, so we can estimate the root in the left-hand side from below to get
	\begin{align*}
		u(a,\theta)
		>
		2 \left( 1 - \cos \theta \frac{\cos \theta + \sqrt{4a^2 - 4a + \cos^2 \theta}}{2a} \right).
	\end{align*}
	Therefore, $u>v$ reduces to
	\begin{align*}
		2
		&>
		\frac{2 a \cos \theta}{\sqrt{4a^2-4a+\cos^2\theta}}
		&& \iff
		\\
		4a^2-4a + \cos^2 \theta
		&>
		a^2 \cos^2 \theta
		&& \iff
		\\
		(a-1)(a(4-\cos^2 \theta)-\cos^2\theta)
		&>
		0
	\end{align*}
	which is obviously true due to~$a>1$.
\end{proof}

Briefly summarizing this section with focus on the information relevant for the main proof in chapter~\ref{chapter proof of main theorem}, we can state the following:

\begin{corollary}\label{cor generating orbits}
	Let~$\gamma$ be an ordinary second species generating orbit with action~$\mathcal{A}(\gamma)$ composed of S-arcs and T-arcs with energy~$c \in (-3/2, \sqrt{2})$ where all supporting ellipses are non-rectilinear and S-arcs have semi-major axes~$a>1$.
	Then there exists an~$\varepsilon>0$ and~$\mu_0 > 0$ such that for all~$\mu \in (0,\mu_0]$ there exists a unique periodic orbit~$\gamma_\mu$ with energy~$c$ in the restricted three-body problem with mass ratio~$\mu$ shadowing~$\gamma$ with action~$\betrag{\mathcal{A}(\gamma_\mu) - \mathcal{A}(\gamma)} < \varepsilon$.
\end{corollary}
Backed by this result, we will in the next chapter compute the action of generating orbits in order to construct orbits with negative action in the restricted three-body problem with positive mass ratio $\mu \in (0,1)$.

%% file: action.tex
In this chapter we compute the action of first and second species generating orbits for the restricted three-body problem.
The first two sections will focus on general formulae for the action of the two species, while the third section will provide a helpful method to compute the elapsed time of second species arcs.
Finally, in the last section we construct sequences of second species generating orbits that have special properties with respect to their action and energy.
First of all however, we will recall the formula for the action and convert to fixed coordinates~$(Q,P)$ in order to subsequently use the geometry of Kepler orbits.

Let $\gamma \colon S^1 \to \Sigma_c$ be a nonconstant periodic orbit of the planar restricted three-body problem with period $\tau > 0$ and energy~$H_\mu = c$.
In rotating coordinates
\begin{align*}
	\begin{pmatrix}
		q_1 \\ q_2 \\ q_3
	\end{pmatrix}
	= 
	\begin{pmatrix}
		\cos(t) & \sin(t) & 0 \\
		- \sin(t) & \cos(t) & 0 \\
		0 & 0 & 1
	\end{pmatrix}
	\begin{pmatrix}
		Q_1 \\ Q_2 \\ Q_3
	\end{pmatrix}
\end{align*}
we have the action
\begin{align} \label{eq action in rotating coordinates}
\begin{split}
\mathcal{A}(\gamma)
:&=
\int_{S^1} \gamma^* \lambda
\\
&=
\int_{0}^{\tau} p_1(t) \frac{\diff q_1(t)}{\diff t} + p_2(t) \frac{\diff q_2(t)}{\diff t} \diff t
\\
&=
\int_{0}^{\tau} p_1(t) \left( p_1(t) + q_2(t) \right) + p_2(t) \left( p_2(t) - q_1(t) \right) \diff t
\\
&=
\int_{0}^{\tau} \norm{p(t)}^2 + L(t) \diff t,
\end{split}
\end{align}
where $L(t)$ is the angular momentum $p_1(t)q_2(t) - p_2(t) q_1(t)$.
Using fixed coordinates $(Q,P)$ allows us to use all our knowledge about Kepler orbits for our generating orbits.
Since both~$\norm{p}^2$ and~$L(t)$ are invariant under the rotation that defines the transformation between fixed and rotating coordinates, the action can simply be written as
\begin{align*}
\mathcal{A}(\gamma)
&= \int_{0}^{\tau} \norm{P(t)}^2 + L(t) \diff t.
\end{align*}
Solving the Kepler Hamiltonian~\eqref{eq Kepler Hamiltonian} for~$\norm{P}^2$ and inserting, we are left with
\begin{align} \label{eq action for general orbits}
\mathcal{A}(\gamma) = \int_{0}^{\tau} 2 H_\mathrm{fix}(t) + L(t) + \frac{2}{\norm{Q(t)}} \diff t.
\end{align}
While for general orbits of the restricted three-body problem the Kepler energy~$H_\mathrm{fix}$, angular momentum~$L$ and distance to the origin~$\norm{Q}$ all depend on time, at least~$H_\mathrm{fix}$ and~$L$ are integrals of motion for the limit case~$\mu=0$.
Since generating orbits are merely rotating Kepler ellipses or sequences of Keplerian arcs, we can now quite easily compute their action.
To further get rid of the term $2/\norm{Q}$, we apply the change of variables to the regularised time~\eqref{eq Levi-Civita time transformation} of the Levi-Civita transformation.
Ultimately, we get that the action of a generating orbit~$\gamma$ which consists of only one Keplerian arc can be computed by
\begin{align} \label{eq action for arcs}
\begin{split}
\mathcal{A}(\gamma)
&= \tau (2 H_\mathrm{fix} + L) + \int_0^\sigma \frac{2}{\norm{X}^2} 4 X^2 \diff s
\\
&= \tau (2 H_\mathrm{fix} + L) + 8 \sigma,
\end{split}
\end{align}
where~$\sigma = s(\tau)$ is the duration of the arc in Levi-Civita-regularised time.
For generating orbits composed of multiple arcs both~$H_\mathrm{fix}$ and~$L$ can change between arcs while the relevant energy~$H_0 = H_\mathrm{fix} + L$ must remain the same.
The total action can simply be computed by adding up the actions of all arcs.

\subsection{First species}
We first compute the action for generating orbits of the fist species, \ie when the orbit is just a full rotating Kepler orbit.
One has to keep in mind, however, that the notion of periodicity remains that from the rotating setting.

\subsubsection{First kind}
For the fist kind---circular orbits---the computation is quite easy.
One only has to compute the period, which in general differs from the Kepler period.

For circular orbits, additionally to the conserved quantities~$H_\mathrm{fix}$ and~$L$, also the radius~$\norm{Q}$ is conserved.
This means that the integrand of~\eqref{eq action for general orbits} itself is a conserved quantity and integration is merely a multiplication with the period:
\begin{align} \label{circular orbits action}
	\mathcal{A}(\gamma) &= \tau \left( 2 H_\mathrm{fix} + L + \frac{2}{\norm{Q}} \right)
\end{align}
Using the computation of the angular velocity~$n := 2 \pi \epsilon' / T = \epsilon' /\sqrt{a^3}$ in fixed coordinates, the mean motion is decreased by the rotation of the coordinate axes to give an angular velocity of $n-1 = \epsilon' /\sqrt{a^3} - 1$.
The period of circular orbits in the rotating frame is hence
\begin{align} \label{eq period of circular orbits}
	\tau = \betrag{ \frac{2 \pi}{\frac{\epsilon'}{\sqrt{a^3}} - 1}} = \frac{2 \pi }{ \betrag{ \frac{1}{\sqrt{a^3}} - \epsilon'}}.
\end{align}
Using equations~\eqref{eq Kepler energy in a} and~\eqref{eq angular momentum in a and epsilon'} for the Kepler energy and angular momentum from the data of the ellipse, we get the action of circular orbits in rotating coordinates as a function of the semi-major axis~$a$ and direction of rotation~$\epsilon'$:
\begin{align}\label{eq circular orbits action a}
\begin{split}
\mathcal{A}(\gamma)
&= \frac{2 \pi}{\betrag{ \frac{1}{\sqrt{a^3}} - \epsilon'}} \left( 2\left( -\frac{1}{2a} \right) - \epsilon' \sqrt{a} + \frac{2}{a} \right)
\\
&= \frac{2 \pi}{\betrag{ \frac{1}{\sqrt{a^3}} - \epsilon'}} \left( \frac{1}{a} - \epsilon' \sqrt{a} \right)
\end{split}
\end{align}
The exceptional case, where the period \eqref{eq period of circular orbits} is undefined, is when $n=1$, \ie $a=1$ and $\epsilon' = +1$.
In this case solutions are stationary in rotating coordinates and lie on the unit circle with the free parameter $\phi_0$.

In general, the action of first kind generating orbits is positive if $\epsilon' = -1$, \ie for all retrograde circular orbits~$I_r$.
Direct orbits exist for energies $H_0 = - 1/ (2a) - \epsilon' \sqrt{a} < -3/2$ and have negative action for radii~$a>1$---direct exterior circular orbits~$I_{de}$---and positive action for~$a<1$---direct interior circular orbits~$I_{di}$---with the action converging towards $-2 \pi$ and $+2 \pi$ at the singularity, respectively.
These orbits of negative action are not that interesting for us, because they only exist in the unbounded component of the Hill's region for energies below the first critical energy.

\subsubsection{Second kind}
First species orbits of the second kind are defined by mutually prime numbers of revolution $I$ of $M_2$ and $J$ of $M_3$ around $M_1$ in fixed coordinates.
By Kepler's third law \eqref{eq Kepler's third law} the period is
$$
\tau = 2 \pi I = J T = 2 \pi \sqrt{a^3} J,
$$
and, hence, the semi-major axis $a$ can be expressed as
\begin{align*}
a = \sqrt[3]{\frac{I^2}{J^2}}.
\end{align*}
By the frequency~\eqref{eq regularised angular frequency} of the regularised ellipse, we know that the period of a full Kepler orbit in Levi-Civita regularised time is $S:= s(T) =\pi / \varpi = \pi /\sqrt{-8H_\mathrm{fix}}$.
The action from \eqref{eq action for arcs} then becomes
\begin{align*}
	\mathcal{A}
	&= J T (2 H_\mathrm{fix} + L) + 8 J S
	\\
	&= 2 \pi \sqrt{a^3} J \left( 2 \left( - \frac{1}{2a} \right) - \epsilon' \sqrt{a(1-\epsilon^2)} \right) + \frac{8 \pi J}{\sqrt{-8 \left( -\frac{1}{2a} \right)}}
	\\
	&= 2 \pi J \left(\sqrt{a} - \epsilon' a^2 \sqrt{1-\epsilon^2} \right)
	\numberthis \label{eq second kind action with a}
	\\
	&= 2 \pi \left(\sqrt[3]{I J^2} - \epsilon' \sqrt[3]{\frac{I^4}{J}} \sqrt{1-\epsilon^2} \right)
\end{align*}
and we see that while there are many combinations of $I$, $J$ and $\epsilon$ that produce negative action, line \eqref{eq second kind action with a} shows that the semi-major axis $a$ must be greater than 1 on order for the action to become negative.
So, in the bounded component~$\mathfrak{H}_c^b$ of the Hill's region no generating orbit of the first species can produce negative action.

\subsection{Second species}
A periodic second species generating orbit is a periodic sequence of Keplerian arcs joined at collision with~$M_2$.
So in order to compute its action we need to compute the action of theses arcs and add them up.
Keplerian arcs are divided into four types as in section \ref{section second species}: Type~1 intersects the unit circle in two distinct points, type~2 is tangent to the unit circle at one point and types~3 and~4 are identical to the unit circle in retrograde and direct direction, respectively.
Type~1 is subdivided into S-arcs and T-arcs, where in S-arcs the collisions occur in two distinct points on the unit circle, while  T-arcs begin and end in collision on the same point on the unit circle in fixed coordinates.
T-arcs can be computed as in the last section, since both $M_2$ and $M_3$ revolve an integer number of times during one arc.
The same holds for arcs of type~2.

For S-arcs we shift time and rotate, such that at $t=0$ the orbit is at it's apocentre or pericentre and on the positive abscissa.
We are going to compute the first intersection time of the regularised orbit with the unit circle in order to find the regularised period of the generating orbit.
Setting the apocentre on the positive abscissa here corresponds to outgoing S-arcs, while the pericentre yields ingoing S-arcs.
Remember that type~1 arcs require the pericentre to be closer than 1 and the apocentre to be farther away than 1 in order to get two distinct intersection points of the Kepler ellipse with the unit circle.
The parametrisation of the regularised orbit can then be given from~\eqref{eq parametrization of regularised Kepler orbit}, \eqref{eq alpha in a and epsilon} and~\eqref{eq beta in a and epsilon} by
\begin{align*}
	X_1(s) = \alpha \cos ( \varpi s)
	\quad \text{and} \quad
	X_2(s) = \beta \sin ( \varpi s),
\end{align*}
where $\alpha = \sqrt{Q_1(0)}$ and $\beta = 2 \dot{Q}_2(0) \sqrt{Q_1(0)} / \varpi$.
The regularised orbit first intersects the unit circle at time $\sigma_0/2 > 0$ when
\begin{align*}
	1 = X_1\left( \frac{\sigma_0}{2} \right)^2 + X_2\left( \frac{\sigma_0}{2} \right)^2.
\end{align*}
Note that $\sigma_0$ is not necessarily the actual regularised period $\sigma$, since an S-arc can first wind around $M_1$ at the origin $J$ times before colliding with $M_2$ again.
The actual regularised period will then be $\sigma = J S + \sigma_0$.
Obviously, the Levi-Civita transformation \eqref{eq Levi-Civita transformation} preservers the unit circle, so an intersection with the unit circle in regularised coordinates corresponds to an intersection in usual coordinates.
Inserting the parametrisation in regularised coordinates, one gets
\begin{align*}
	&&1 &= \alpha^2 \cos^2 \left( \varpi \frac{\sigma_0}{2} \right) + \beta^2 \sin^2 \left( \varpi \frac{\sigma_0}{2} \right)
	\\
	&&&= \alpha^2 \cos^2 \left( \varpi \frac{\sigma_0}{2} \right) + \beta^2 \sin^2 \left( \varpi \frac{\sigma_0}{2} \right) + \beta^2 \cos^2 \left( \varpi \frac{\sigma_0}{2} \right) - \beta^2 \cos^2 \left( \varpi \frac{\sigma_0}{2} \right)
	\\
	&&&= \beta^2 + \cos^2 \left( \varpi \frac{\sigma_0}{2} \right) (\alpha^2 - \beta^2)
	\\
	\iff&&
	\cos^2 \left( \varpi \frac{\sigma_0}{2} \right)
	&= \frac{1 - \beta^2}{\alpha^2 - \beta^2}
	\\
	\iff&&
	\sigma_0
	&= \frac{2}{\varpi} \arccos \left( \sqrt{\frac{1 - \beta^2}{\alpha^2 - \beta^2}} \right).
\end{align*}
Equivalences hold because $\varpi \sigma_0 / 2 \in (0,pi)$ and $\beta^2 = \alpha^2 \iff \alpha = 1 \land \beta = \pm 1$.
Both these latter cases do not admit an S-arc, but are rather of types~3 and~4, respectively.
Using the formulae~\eqref{eq regularised angular frequency} for the angular frequency, \eqref{eq alpha in a and epsilon} for~$\alpha$, and \eqref{eq beta in a and epsilon} for~$\beta$ in terms of the Kepler ellipse data, we get
\begin{align*}
	\sigma_0
	&= \frac{2}{\sqrt{-8H_\mathrm{fix}}} \arccos \left(\sqrt{\frac{1 - \beta^2}{\alpha^2 - \beta^2}}\right)
	\\
	&= \frac{2}{\sqrt{-8\left(- \frac{1}{2a}\right)}} \arccos \left(\sqrt{\frac{1 - a(1 - \epsilon)}{a(1 + \epsilon) - a(1 - \epsilon)}}\right)
	\\
	&= \sqrt{a} \arccos \left(\sqrt{\frac{1 - a(1 - \epsilon)}{2 a \epsilon}}\right).
	\numberthis \label{eq sigma_0 with arccos}
\end{align*}
\begin{figure}[tb]
\centering
\includegraphics{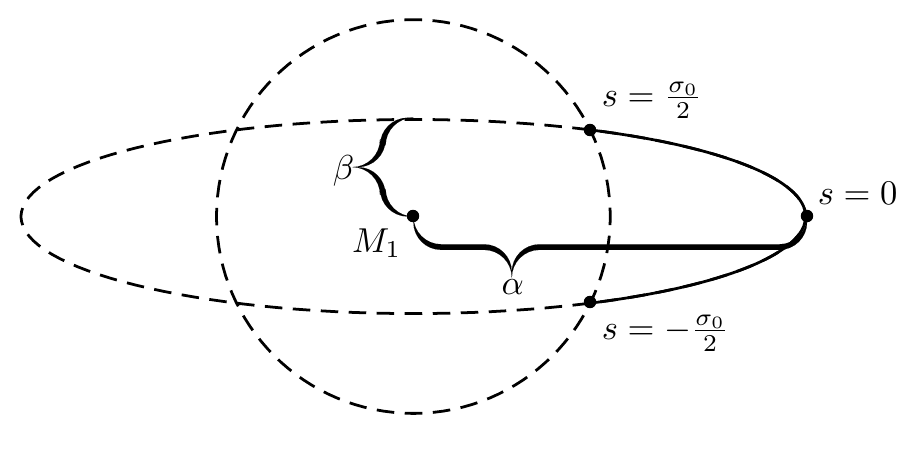}
\caption{First intersection time of Levi-Civita-regularised Kepler ellipses with the unit circle.}
\end{figure}
Inserting this into the formula \eqref{eq action for arcs} for the action, we get
\begin{flalign*}
	\mathcal{A}
	&= (J T + \tau_0) (2H+L) + 8 (J S + \sigma_0)
	\\
	&= (2 \pi \sqrt{a} J + \tau_0) \left( 2 \left( - \frac{1}{2a} \right) - \epsilon' \sqrt{a(1-\epsilon^2)} \right) + 8 \left( J \frac{\pi}{\varpi} + \sqrt{a} \arccos \left( \sqrt{\frac{1 - a(1 - \epsilon)}{2 a \epsilon}} \right) \right) 
	\\
	&= (2 \pi \sqrt{a} J + \tau_0) \left(- \frac{1}{a} - \epsilon' \sqrt{a(1-\epsilon^2)} \right) + 8 \sqrt{a} \left( \frac{\pi J}{2} + \arccos \left( \sqrt{\frac{1 - a(1 - \epsilon)}{2 a \epsilon}} \right) \right),
	\numberthis \label{eq action of second species}
\end{flalign*}
where $\tau_0 = t(\sigma_0)$ is two times the first intersection time with the unit circle in normal time.
This quantity has to be computed separately, which is done in the subsequent section.

The action of type~2 arcs is identical to the previous case of first species generating orbits of the second kind, since they are complete integer revolutions of Keplerian ellipses.
Type~3 is half a circular retrograde Kepler orbit with radius $a=1$ and identical with the circular retrograde generating orbit of family $I_{r}$ with radius 1 and action $\mathcal{A} = 2 \pi$.
Type~4 on the other hand is in rotating coordinates the constant solution identical with $M_2$ at all times and therefore a third species generating orbit, which can not be computed using the Kepler problem.

\subsection{Lambert's Theorem} \label{section Lamberts theorem}
A great tool that especially enables us to compute the elapsed time of Keplerian arcs of the second species is Lambert's theorem from~\cite{Lambert}.
Its history, modern proofs and many remarks about it can be found in \cite{Albouy_Lambert} and \cite{AlbouyUrena_simple}.
We will state the main theorems and definitions needed for this work, while adapting the notation slightly.
Objects of study for Lambert's theorem are Keplerian arcs beginning in point~$A$ at time~$t_A$ and ending in point~$B$ at time~$t_B$ in fixed coordinates.

\begin{theorem}[Lambert, Theorem 1 in \cite{Albouy_Lambert}] \label{theorem Lambert 1}
	Consider Keplerian arcs around the origin $O$ of $\RR^d$.
	If we change continuously such an arc while keeping constant the distance~$\norm{B-A}$ between both ends, the sum of the radii~$\norm{A} + \norm{B}$ and the energy~$H_\mathrm{fix}$, then the elapsed time $\tau_0 = t_B - t_A$ is also constant.
\end{theorem}
\begin{theorem}[Lambert, Theorem 2 in \cite{Albouy_Lambert}] \label{theorem Lambert 2}
	Starting from any given Keplerian arc, we can arrive at some rectilinear arc by a continuous change which keeps constant the three quantities $\norm{B-A}$, $\norm{A} + \norm{B}$ and~$H_\mathrm{fix}$.
\end{theorem}
\begin{figure}[tb]
\centering
\includegraphics{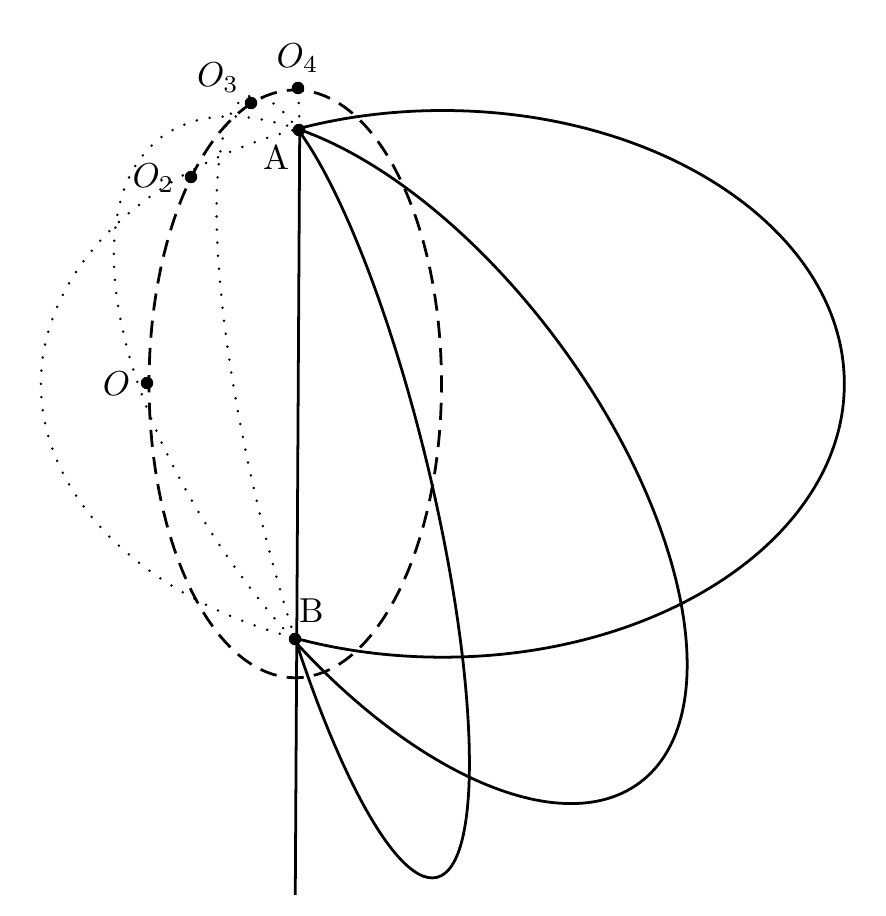}
\caption{Lambert cycle of Keplerian arcs where the origin moves along a second ellipse with foci~$A$ and~$B$.}
\end{figure}
\begin{definition}[see Definition 5 in \cite{Albouy_Lambert}] \label{definition Lambert direct indirect}
	A Keplerian arc around~$O$ is called \textit{simple} if its elapsed time is less than or equal to the period of its supporting ellipse.
	It is said to be \textit{indirect}, or~$I_O$, if its convex hull contains~$O$;
	\textit{direct}, or~$D_O$, if its convex hull does not contain~$O$;
	\textit{indirect with respect to the second focus~$F$}, or~$I_F$, if its convex hull contains~$F$;
	\textit{direct with respect to~$F$}, or~$D_F$, if its convex hull does not contain~$F$.
\end{definition}
In order to avoid confusion between terminologies we will only use the notation of $I_O$, $D_O$, $I_F$ and $D_F$.
The important feature of these types is that they are preserved during the \textit{Lambert cycle}, which is the continuous change of Keplerian arcs as in theorem \ref{theorem Lambert 1}:
\begin{proposition}[Proposition 5 in \cite{Albouy_Lambert}] \label{proposition Lambert}
	If a Keplerian arc in a Lambert cycle is~$I_O$ or~$D_O$ and~$I_F$ or~$D_F$, then all the Keplerian arcs of the cycle are~$I_O$ or~$D_O$ and~$I_F$ or~$D_F$, respectively.
\end{proposition}

In our case of second species generating orbits we have as parameters the semi-major axis~$a$, the eccentricity~$\epsilon$ and the polar angle~$\theta$ of the intersection with the unit circle.
They are interrelated by the equation for the polar distance

\begin{align} \label{eq focal distance of ellipse}
	r(\theta) = \frac{a(1-\epsilon^2)}{1-\epsilon \cos \theta}
\end{align}
of an ellipse with one focus in the origin, so we have

\begin{align*}
\cos \theta = \frac{1-a(1-\epsilon^2)}{\epsilon}.
\end{align*}
The elapsed time can then be computed with the help of theorem~\ref{theorem Lambert 1} and~\ref{theorem Lambert 2} by the elapsed time of a rectilinear arc.
Since the energy~$H_\mathrm{fix} = - 1/(2a)$ remains constant, the semi-major axis remains constant during the continuous change and, consequently, the supporting rectilinear orbit has length~$2a$.
We can follow the shape of the supporting ellipse during the continuous change by shifting the focus of the ellipse along a second ellipse with foci~$A$ and~$B$ going through the origin.
The semi-major axis of this second ellipse is 1 in our case, since~$A$ and~$B$ lie on the unit circle, and therefore the elapsed time can be computed by using the free-fall time~\eqref{eq free-fall time} from height~$2a$ to~$A$ and to~$B$.

By Proposition~\ref{proposition Lambert}, depending on whether the original arc was~$I_O$ or~$D_O$ and~$I_F$ or~$D_F$, the new rectilinear arc will also be~$I_O$ or~$D_O$ and~$I_F$ or~$D_F$, respectively.
The original arc is
\begin{align*}
\left\lbrace
\begin{aligned}
I_O \quad &\text{if }\theta \geq \frac{\pi}{2}
\\
D_O \quad &\text{if }\theta < \frac{\pi}{2}
\end{aligned}
\right\rbrace
\qquad \text{and} \qquad
\left\lbrace
\begin{aligned}
I_F \quad &\text{if } 2 a \epsilon \geq \cos \theta
\\
D_F \quad &\text{if } 2 a \epsilon < \cos \theta
\end{aligned}
\right\rbrace,
\end{align*}
and we can state the following conclusion:
\begin{lemma} \label{lemma time of arc}
The elapsed time of an outgoing second species generating arc is
\begin{align*}
	\tau
	&= 2 \pi J \sqrt{a^3} +
	\begin{cases}
		t_{2a} ( 1 + \sin \theta )
		+ t_{2a}(0)
		+ ( t_{2a}(0) - t_{2a} ( 1 - \sin \theta ) )
		\quad &\text{if } \cos \theta \leq 0
		\\
		t_{2a} ( 1 + \sin \theta ) + t_{2a} ( 1 - \sin \theta )
		\quad &\text{if } 0 < \cos \theta \leq 2 a \epsilon
		\\
		t_{2a} ( 1 + \sin \theta ) - t_{2a} ( 1 - \sin \theta )
		\quad &\text{if } 2 a \epsilon < \cos \theta.
	\end{cases}
\end{align*}
\end{lemma}

\begin{remark}
	There is no case~$I_O$ and~$D_F$.
	If we rotate the arc such that the apoapsis---which has to be farther away than 1---lies on the $Q_1$-axis, then the second focus lies between the apoapsis and the origin on the $Q_1$-axis.
	Hence, if an outgoing arc is~$I_O$, it is necessarily also~$I_F$.
	The times of ingoing second species arcs are simply $2\pi \sqrt{a^3} = 2 t_{2a}(0)$ minus the outgoing time.
\end{remark}
The orbits which we will construct in the next section will only be outgoing second species generating orbits without~$M_3$ winding around~$M_1$, \ie with $J=0$.
This restriction makes sense particularly in view of~\eqref{eq action for general orbits}, where for negative angular momentum the only positive contribution to the action comes from the the term $1/\norm{Q}$ which we try to keep as small as possible in order to get negative action.
In this situation the following statement helps us to exclude orbits without winding of~$M_2$ around~$M_1$.
\begin{theorem}[Theorem 7.2 from \cite{AlbouyUrena_simple}] \label{theorem uniqueness of arcs}
	In the Euclidean plane or space consider three distinct points~$O$, $A$, $B$ such that~$O$ is not on the segment $AB$.
	There is a unique $D_O$ Keplerian arc around~$O$ and a unique simple $I_O$ Keplerian arc around~$O$ going from A to B in a given positive elapsed time.
	In the exceptional case $O \in ]A,B[$ there exist exactly two distinct $I_O$ Keplerian arcs which are reflections of each-other.
\end{theorem}
\begin{corollary}\label{corollary I not zero}
	There exist no S-arcs with $I=J=0$ and $\epsilon'=+1$ which are not identical to $M_2$ at all times.
\end{corollary}
\begin{proof}
An S-arc requires a Keplerian arc with two distinct ends $A$ and $B$ on the unit circle.
The timing condition for $I=J=0$ requires that the elapsed time of the arc is exactly the elapsed time of~$M_2$ between~$A$ and~$B$.
Since both~$M_2$ and~$M_3$ move in the same direct direction around~$M_1$, the arc of~$M_3$ is $D_O$ and~$I_O$ if and only if the arc of $M_2$ is~$D_O$ and~$I_O$, respectively.
Uniqueness in theorem \ref{theorem uniqueness of arcs} give us that both arcs are the same and hence~$M_3$ coincides with~$M_2$ for all time.
\end{proof}

\subsection{Sequences of generating orbits} \label{section sequences of generating orbits}

We will now describe some sequences of generating orbits and also show where they are found in terms of Hénon's notation for generating families.
Using our formulae from the previous sections of this chapter, we can find generating orbits with negative action that have additional properties.
In our case we want to control the energy and show that these orbits exist for all energies between $-\sqrt{2}$ and 0.
All generating orbits described here will be of the second species, but not all orbits in this section will actually have negative action.
They are then rather included here either because they are instructive and arise in Hénon's classification of families or because they are at the beginning of a sequence where the action tends towards negative values but is not necessarily negative throughout the entire sequence.

The main formula we will use to compute the action is equation~\eqref{eq action of second species}, which is quite powerful but still requires the elapsed time of the Kepler arc between collisions with $M_2$.
Since the Kepler arc corresponding to the generating orbit is usually not a full Kepler ellipse, we will use Lemma~\ref{lemma time of arc} for the remaining cases.
To further simplify things, we will always set $J = 0$, \ie the arc of $M_3$ does not wind around $M_1$ in fixed coordinates.
The reason for this is equation~\eqref{eq action for general orbits}, where the only term contributing positively to the action is $1/\norm{Q}$ and we intend to keep this term small by not letting the orbit come unnecessarily close to~$M_1$.
Also, our generating orbits here will only consist of a single outgoing $S$-arc which is repeated infinitely to give a periodic generating orbit.

Another issue in order to analytically describe second species generating orbits is the timing condition
\begin{align} \label{eq timing condition}
	\tau = \begin{cases}
		2 \pi I + 2 \theta &\qquad \text{for direct and rectilinear orbits}
		\\
		2 \pi (I+1) - 2 \theta &\qquad \text{for retrograde orbits,}
	\end{cases}
\end{align}
\ie the requirement that~$M_3$ collides again with~$M_2$ after starting at collision and following the arc.
This problem is avoided in this chapter by leaving a free parameter which will be the semi-major axis~$a$ of the supporting Kepler ellipse.
We will then show that~$\tau - 2 \epsilon' \theta$ tends smoothly towards infinity as $a \to \infty$.
This means that for all large enough numbers of revolution~$I$ of~$M_2$ around~$M_1$ in fixed coordinates there will be an~$a$ that solves the timing condition for that particular~$I$.
A one-parameter family with~$a$ as the free parameter will in this way give a sequence of second species generating orbits with an orbit for each large enough integer~$I \geq 1$.
Other parameters of the supporting Keplerian orbit which will be used in this section are the eccentricity~$\epsilon$, the semi-minor axis~$b$, the polar angle of intersection with the unit circle~$\theta$ and the direction of rotation~$\epsilon'$.

\subsubsection{Fixed semi-minor axis}
The first sequence we want to present is one where the energy converges towards zero and the action against negative infinity.
This can be achieved by fixing the semi-minor axis~$b$.
The relation between~$b$, $a$ and $\epsilon$ is 
$b = a \sqrt{1-\epsilon^2}$, \ie
\begin{align} \label{eq epsilon fixed b}
	\epsilon = \sqrt{1- \frac{b^2}{a^2}}.
\end{align}
Inserting this into the Kepler energy~\eqref{eq Kepler energy in a} and the angular momentum \eqref{eq angular momentum in a and epsilon'}, the rotating Kepler energy~\eqref{eq rotating Kepler Hamiltonian} becomes
\begin{align}
	H_0 = - \frac{1}{2a} - \epsilon' \frac{b}{\sqrt{a}}
\end{align}
which strictly monotonically tends towards zero from below as $a \to \infty$ and $\epsilon'=+1$ for any fixed~$b$.
We choose the easiest nonzero $b=1$ for the sequence.

In order to define an outgoing arc not part of the unit circle, we need the maximal focal distance
\begin{align*}
	a(1+\epsilon) = a + \sqrt{a^2-1} > 1,
	\qquad \text{\ie} \qquad a > 1.
\end{align*}

For the elapsed time of the arc we use Lemma \ref{lemma time of arc}, for which we need $\sin \theta$.
Also, we need to check if the arcs are~$I_O$ and~$I_F$, or~$D_O$ and~$I_F$, or~$D_O$ and~$D_F$.
We can compute~$\cos \theta$ from the equation for focal distances of ellipses~\eqref{eq focal distance of ellipse}:
\begin{align*}
	\cos \theta
	&=
	\frac{1 - a(1-\epsilon^2)}{\epsilon}
	=
	\frac{1 - \frac{1}{a}}{\sqrt{1- \frac{1}{a^2}}}
	=
	\sqrt{\frac{a-1}{a+1}}
\end{align*}
Since~$a>1$, we have
\begin{align*}
	0
	< \cos \theta
	= \sqrt{\frac{a-1}{a+1}}
	< 2 \sqrt{(a-1)(a+1)}
	= 2 a \epsilon,
\end{align*}
\ie all arcs in this sequence are~$D_O$ and~$I_F$, and we can compute
\begin{align*}
	\tau
	&=t_{2a}(1-\sin \theta) + t_{2a}(1+\sin \theta)
	\\
	&=2 \sqrt{a^3}\begin{aligned}[t]
		&\left(
		\sqrt{\frac{1 - \sin \theta}{2a} \left( 1 - \frac{1 - \sin \theta}{2a} \right) }
		+ \sqrt{\frac{1 + \sin \theta}{2a} \left( 1 - \frac{1 + \sin \theta}{2a} \right) }
		\right. \\
		&\left.
		+ \arccos \left( \sqrt{\frac{1 - \sin \theta}{2a}} \right)
		+ \arccos \left( \sqrt{\frac{1 + \sin \theta}{2a}} \right)
		\right),
	\end{aligned}
	\numberthis \label{eq tau fixed b}
\end{align*}
where
\begin{align} \label{eq sin theta fixed b}
	\sin \theta
	&= \sqrt{1 - \cos^2 \theta}
	= \sqrt{1 - \frac{a-1}{a+1}}
	= \sqrt{\frac{2}{a+1}}.
\end{align}
Obviously, the elapsed time tends towards infinity as~$a \to \infty$.
In particular for~$\epsilon'=+1$ the time~$\tau - 2 \theta$ becomes zero for the limit case~$a=1$ and depends smoothly on~$a$.
So there exists for every $I \geq 1$ an~$a>1$ solving the timing condition~\eqref{eq timing condition} and we get a sequence of second species generating orbits.

The regularised time from~\eqref{eq sigma_0 with arccos} becomes
\begin{align*}
	\sigma
	&= \sqrt{a} \arccos \left( \sqrt{\frac{1 - a(1-\epsilon)}{2a\epsilon}} \right) 
	\\
	&= \sqrt{a} \arccos \left( \sqrt{\frac{1 - a \left( 1-\sqrt{1- \frac{1}{a^2}} \right) }{2\sqrt{a^2-1}}} \right) 
	\\
	&= \sqrt{a} \arccos \left( \sqrt{\frac{1}{2} - \frac{1}{2} \sqrt{\frac{a-1}{a+1}}} \right)
	\numberthis \label{eq sigma fixed b}
\end{align*}
and by inserting~\eqref{eq epsilon fixed b},~\eqref{eq tau fixed b} and~\eqref{eq sigma fixed b} into \eqref{eq action for arcs}, we get the action of these generating orbits:
\begin{align*}
	\mathcal{A}
	&= 2 \sqrt{a^3}\left( - \frac{1}{a} - \epsilon' \frac{1}{\sqrt{a}} \right)
		\left(
		\sqrt{\frac{1 - \sqrt{\frac{2}{a+1}}}{2a} \left( 1 - \frac{1 - \sqrt{\frac{2}{a+1}}}{2a} \right) }
		\right. \\
		&\left.
		+ \sqrt{\frac{1 + \sqrt{\frac{2}{a+1}}}{2a} \left( 1 - \frac{1 + \sqrt{\frac{2}{a+1}}}{2a} \right) }
		+ \arccos \left( \sqrt{\frac{1 - \sqrt{\frac{2}{a+1}}}{2a}} \right)
		+ \arccos \left( \sqrt{\frac{1 + \sqrt{\frac{2}{a+1}}}{2a}} \right)
		\right)
	\\
	& \hphantom{=} + 8 \sqrt{a} \arccos \left( \sqrt{\frac{1}{2} - \frac{1}{2} \sqrt{\frac{a-1}{a+1}}} \right)
\end{align*}
This tends towards negative infinity as~$a$ goes towards infinity.

Summarising so far, we can state
\begin{lemma} \label{lemma sequence zero energy}
	There exists a sequence of second species generating orbits consisting of non-rectilinear S-arcs with semi-major axis $a>1$, action tending towards negative infinity and energy~$H_0$ strictly monotonically converging towards zero from below.
\end{lemma}

\subsubsection{Fixed polar intersection angle with the unit circle}\label{subsection fixed theta}

For more general sequences we fix the angle of intersection~$\theta$ with the unit circle.
We will look at all~$\theta \in [0,\pi]$ but we will also highlight some special cases that arise.
Let from now on~$a>1$.
The eccentricity is computed by solving the equation~\eqref{eq focal distance of ellipse} of the focal distance of ellipses for~$\epsilon$:
\begin{align*}
	1 &= \frac{a(1-\epsilon^2)}{1-\epsilon \cos \theta}
	\\ \iff \quad
	a \epsilon^2 - \epsilon \cos \theta - a + 1 &= 0
	\\ \iff \quad
	\epsilon &= \frac{\cos \theta \pm \sqrt{\cos^2 \theta + 4a(a-1)}}{2a}
	\numberthis \label{eq epsilon fix theta both signs}
	\\ \overset{\epsilon \geq 0}{\iff} \quad
	\epsilon &= \frac{\cos \theta + \sqrt{\cos^2 \theta + 4a(a-1)}}{2a}
	\numberthis \label{eq epsilon fix theta}
\end{align*}
\begin{remark} \label{remark fixed theta a greater 1}
	There exist ellipses intersecting the unit circle at polar angle~$\theta \in [0,\pi/2]$ also for semi-major axes~$a \leq 1$.
	However,~$a$ is no longer a monotone parameter there.
	Now both signs of~\eqref{eq epsilon fix theta both signs} return nonnegative eccentricity.
	The actually smallest semi-major axis is attained at the largest root of the discriminant~$\cos^2 \theta + 4a(a-1)$.
	A better parameter for small~$a$ would for example be the eccentricity~$\epsilon$---see also proposition 2.5 in \cite{AlbouyUrena_simple}---but we will stick to the semi-major axis as our parameter because computations of the energy and action are easier.
	Another advantage of setting~$a>1$ is that
	\begin{align} \label{eq fixed theta are I_F}
		2a \epsilon = \cos \theta + \sqrt{\cos^2 \theta + 4a(a-1)} > \cos \theta,
	\end{align}
	\ie all arcs are~$I_F$.
\end{remark}
\begin{figure}[tb]
\centering
\begin{subfigure}[t]{0.495\textwidth}
	\centering
	\includegraphics[width=\textwidth]{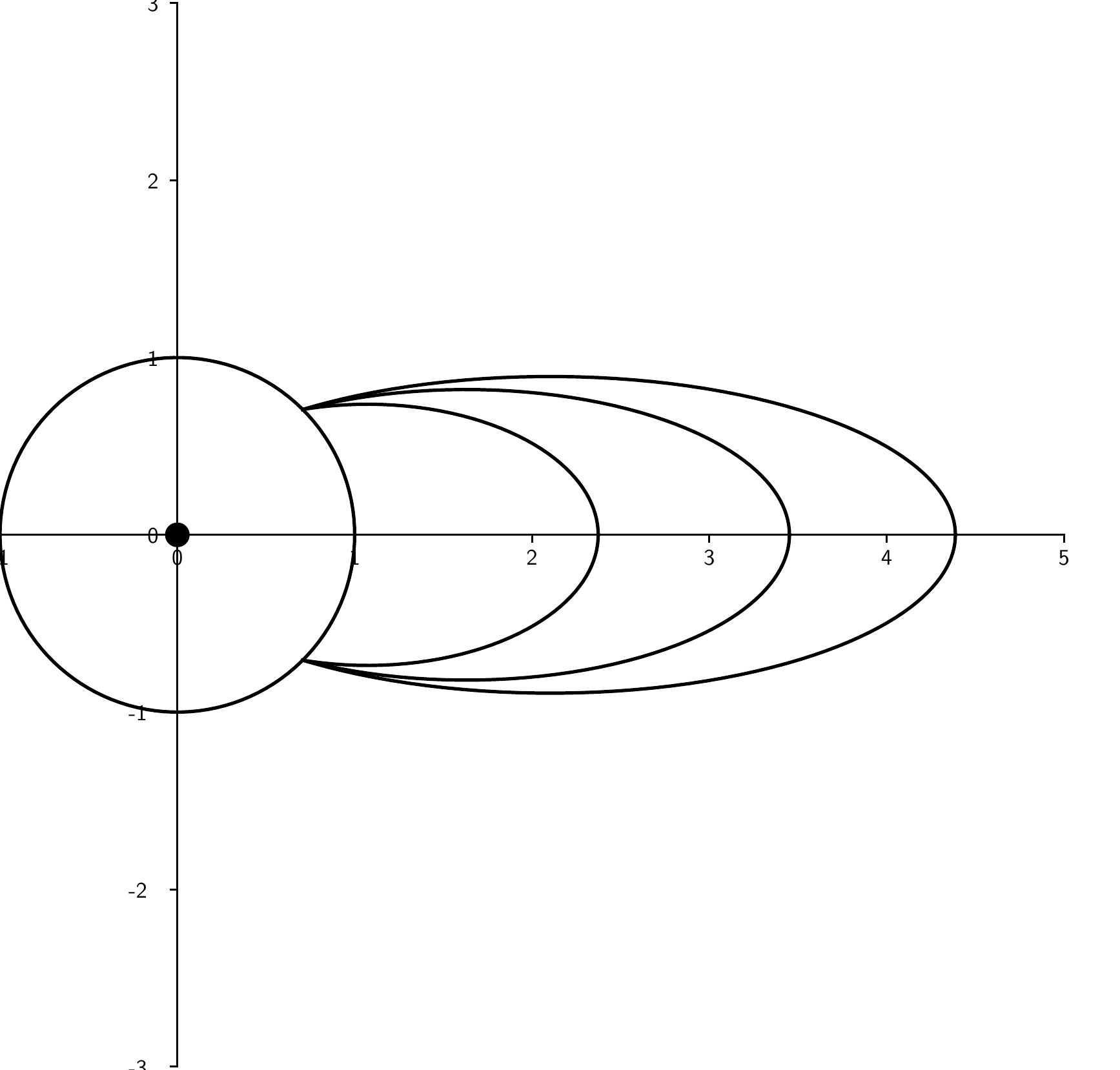}
	\caption{In fixed coordinates.}
\end{subfigure}
\hfill
\begin{subfigure}[t]{0.495\textwidth}
	\centering
	\includegraphics[width=\textwidth]{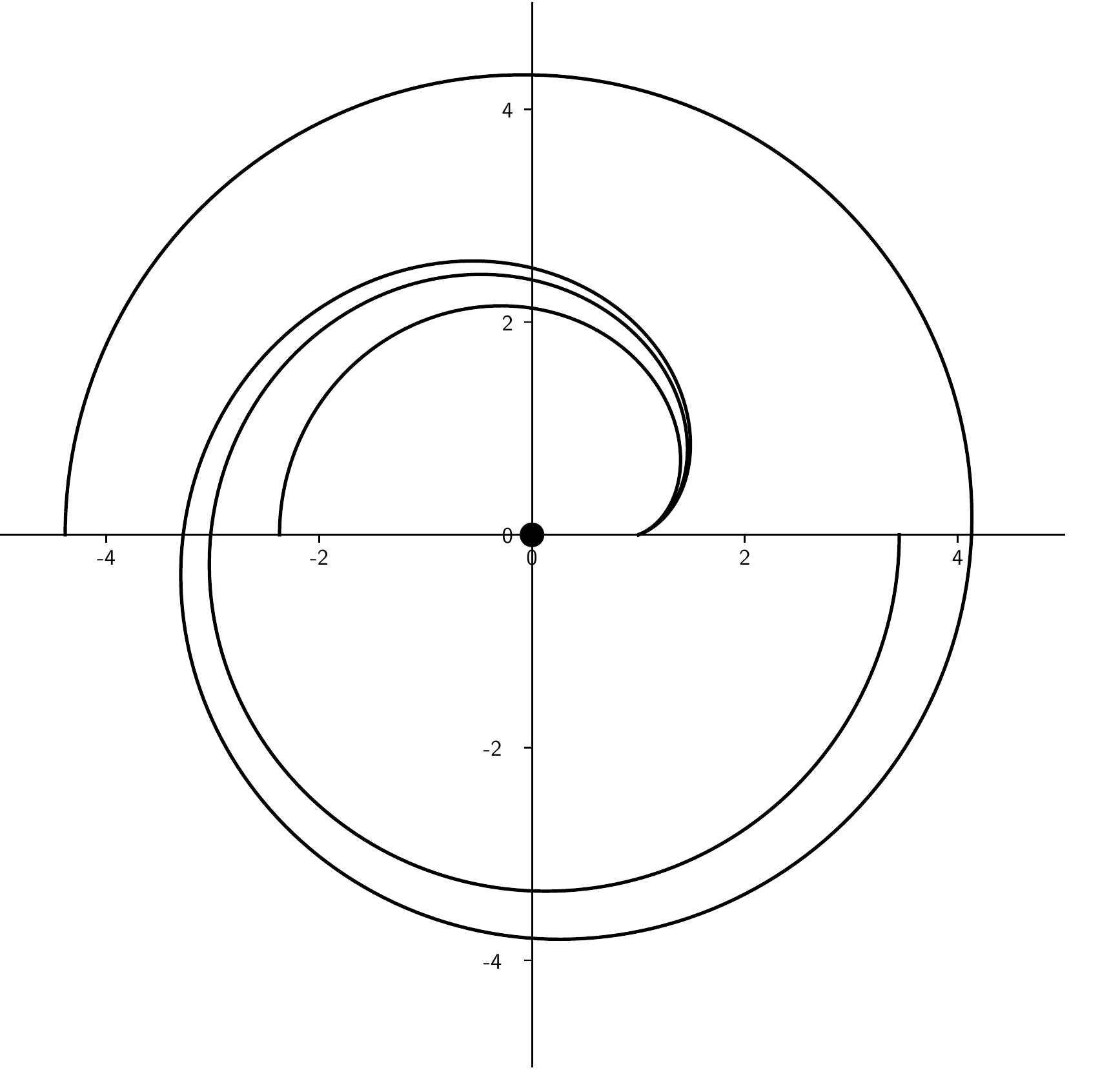}
	\caption{Half orbits in rotating coordinates.}
\end{subfigure}
\caption{Second species arcs with fixed~$\theta = \pi/4$ and~$I = 1,2$.}
\end{figure}
The energy becomes
\begin{align*}
	H_0 &= - \frac{1}{2a} - \epsilon' \sqrt{a(1-\epsilon^2)}
	\\
	&=  - \frac{1}{2a} - \epsilon' \sqrt{a \left( 1 - \frac{\left(\cos \theta + \sqrt{\cos^2 \theta + 4a(a-1)}\right)^2}{4a^2} \right)}
	\\
	&=  - \frac{1}{2a} - \epsilon' \sqrt{\frac{4a^2 - \cos^2 \theta - 2 \cos \theta \sqrt{\cos^2 \theta + 4a(a-1)} - \cos^2 \theta - 4a(a-1)}{4a}}
	\\
	&=  - \frac{1}{2a} - \epsilon' \sqrt{ 1 - \cos \theta  \frac{\cos \theta + \sqrt{4a^2 - 4a  + \cos^2 \theta}}{2a}}
	\quad \overset{a \to \infty}{\longrightarrow} \quad
	- \epsilon' \sqrt{1- \cos \theta}
	\numberthis \label{eq energy fixed theta}
\end{align*}
and tends towards values between~$-\sqrt{2}$ and~$\sqrt{2}$ depending on~$\theta$ and~$\epsilon'$.
It was also computed in lemma~\ref{lemma non-degenerate S-arcs} that the derivative $\partial H_0 / \partial a$ never vanishes, so the energy strictly increases for~$\epsilon' = +1$ and strictly decreases for~$\epsilon' = -1$ along~$a$.

In the first case~$\theta < \pi/2$ all arcs are~$D_O$ and~$I_F$, and we compute the elapsed time from lemma~\ref{lemma time of arc} to be
\begin{align*}
	\tau
	&=
	t_{2a} ( 1 + \sin \theta ) + t_{2a} ( 1 - \sin \theta )
	\\
	&=2 \sqrt{a^3}\begin{aligned}[t]
		&\left(
		\sqrt{\frac{1 - \sin \theta}{2a} \left( 1 - \frac{1 - \sin \theta}{2a} \right) }
		+ \sqrt{\frac{1 + \sin \theta}{2a} \left( 1 - \frac{1 + \sin \theta}{2a} \right) }
		\right. \\
		&\left.
		+ \arccos \left( \sqrt{\frac{1 - \sin \theta}{2a}} \right)
		+ \arccos \left( \sqrt{\frac{1 + \sin \theta}{2a}} \right)
		\right).
	\end{aligned}
	\numberthis \label{eq tau fixed small theta}
\end{align*}
In the second case~$\theta \geq \pi/2$ all arcs are~$I_O$ and~$I_F$, and we use the first case in lemma~\ref{lemma time of arc}, so
\begin{align*}
	\tau
	&= t_{2a} ( 1 + \sin \theta )
	+ t_{2a}(0)
	+ ( t_{2a}(0) - t_{2a} ( 1 - \sin \theta ) )
	\\
	&=2 \sqrt{a^3}
	\begin{aligned}[t]
		&\left(
		\sqrt{\frac{1 + \sin \theta}{2a} \left( 1 - \frac{1 + \sin \theta}{2a} \right) }
		+ \arccos \left( \sqrt{\frac{1 + \sin \theta}{2a}} \right)
		+ \frac{\pi}{2}
		\right. \\
		&\left.
		+ \left( \frac{\pi}{2}
		- \sqrt{\frac{1 - \sin \theta}{2a} \left( 1 - \frac{1 - \sin \theta}{2a} \right) }
		- \arccos \left( \sqrt{\frac{1 - \sin \theta}{2a}} \right)
		\right)
		\right).
	\end{aligned}
\end{align*}
It is obvious in both cases that for fixed~$\theta$ the elapsed time~$\tau$ increases strictly and tends towards infinity as~$a \to \infty$ because~$\tau$ is only the sum of free-fall times to the same points from increased heights~$2a$.
This means the timing condition~\eqref{eq timing condition} has a unique solution for every large enough~$I\geq 1$.
Note that corollary~\ref{corollary I not zero} states that there is no nontrivial solution for~$I=0$.
We also see that the assumption~$a>1$ is no restriction, since for $a=1$ and~$\theta>0$ we get
\begin{align*}
	\tau < 2t_2(0) = 2\pi < 2 \pi + 2 \theta.
\end{align*}
Hence, the timing condition for~$\epsilon' = +1$, $I=1$ and any fixed~$\theta>0$, is always satisfied for an~$a>1$.
The case~$\theta=0$ will later be treated separately.

The regularised time is
\begin{align*}
	\sigma
	&= \sqrt{a} \arccos \left( \sqrt{\frac{1 - a(1-\epsilon)}{2a\epsilon}} \right)
	\\
	&= \sqrt{a} \arccos \left( \sqrt{\frac{\frac{2 - 2a + \cos \theta + \sqrt{\cos^2 \theta + 4a(a-1)}}{2}}{\cos \theta + \sqrt{\cos^2 \theta + 4a(a-1)}}} \right)
	\\
	&= \sqrt{a} \arccos \left( \sqrt{\frac{1}{2} - \frac{a-1}{\cos \theta + \sqrt{\cos^2 \theta + 4a(a-1)}}} \right)
\end{align*}
and the action from~\eqref{eq action of second species} becomes
\begin{align} \label{eq action fixed theta}
	\begin{split}
		\mathcal{A}
		&= \tau \left( -\frac{1}{a} - \epsilon' \sqrt{ 1 - \cos \theta  \frac{\cos \theta + \sqrt{4a^2 - 4a  + \cos^2 \theta}}{2a}} \right)
		\\
		&\hphantom{=} + 8 \sqrt{a} \arccos \left( \sqrt{\frac{1}{2} - \frac{a-1}{\cos \theta + \sqrt{\cos^2 \theta + 4a(a-1)}}} \right).
	\end{split}
\end{align}
The important feature is that the action tends towards negative infinity for every fixed and positive~$\theta$ as~$a \to \infty$ and~$\epsilon'=+1$ since the angular momentum tends towards~$\sqrt{1-\cos \theta}$.
This holds in both cases~$\theta<\pi/2$ and~$\theta \geq \pi/2$ since~$\sqrt{a^3}$ multiplied by bounded terms outweighs the only positive term of~$\sqrt{a}$ times something bounded.

In the limit case~$\theta = 0$, however, the angular momentum vanishes completely and we get
\begin{align*}
	H_0 &= - \frac{1}{2a}
	\\
	\tau &= 2 \sqrt{a^3}
	\left(
	2 \sqrt{\frac{1}{2a} \left( 1 - \frac{1}{2a} \right) }
	+ 2 \arccos \left( \sqrt{\frac{1}{2a}} \right)
	\right)
	\\
	\sigma &= \sqrt{a} \arccos \left( \sqrt{\frac{1}{2a}} \right)
\end{align*}
and hence
\begin{align*}
	\mathcal{A} &=- \frac{4}{a} \sqrt{a^3}
	\left(
	\sqrt{\frac{1}{2a} \left( 1 - \frac{1}{2a} \right) }
	+ \arccos \left( \sqrt{\frac{1}{2a}} \right)
	\right) + 8 \sqrt{a} \arccos \left( \sqrt{\frac{1}{2a}} \right)
	\\
	&= 4 \sqrt{a} \arccos \left( \sqrt{\frac{1}{2a}} \right)
	- 2\sqrt{2 - \frac{1}{a} }.
\end{align*}
The action here no longer tends towards negative values and is in fact always nonnegative.
In fixed coordinates~$M_3$ would fall freely towards~$M_1$---where we would have to regularise---and then bounce back to the place it started.
In our setting of an outgoing generating orbit of the second species with $J = 0$, however, $M_3$ would not make it that far since it would first collide with~$M_2$ moving on the unit circle.

Two more special cases will be mentioned here:
the case~$\theta = \pi/2$ and~$\theta=\pi$.
For~$\theta=\pi$ the arcs are full Kepler ellipses, \ie second species of type 2.
We get
\begin{align*}
	H &= - \frac{1}{2a} - \epsilon' \sqrt{2 - \frac{1}{a}}
	\\
	\tau &= 2 \pi \sqrt{a^3}
	\\
	\sigma &= \frac{\pi \sqrt{a}}{2}
\end{align*}
and
\begin{align*}
	\mathcal{A}
	&=
	2 \pi \sqrt{a^3} \left( - \frac{1}{a} - \epsilon' \sqrt{2 - \frac{1}{a}} \right) + 8 \frac{\pi \sqrt{a}}{2}
	\\
	&= 2 \pi \sqrt{a} - 2 \pi \epsilon' a \sqrt{2a - 1}.
\end{align*}
This action obviously has the same sign as~$-\epsilon'$ for all~$a>1$.
Since these generating orbits are both first and second species orbits, they are naturally bifurcation orbits.
In the continuation to the restricted three-body problem the intersection of the corresponding families splits into two separate families here.

In the remaining special case we fix~$\theta = \pi/2$.
The data simplifies to
\begin{align*}
	H_0 &= - \frac{1}{2a} - \epsilon',
	\\
	\tau &= 2 \sqrt{a^3}
	\left(
	\sqrt{\frac{1}{a}\left(1 - \frac{1}{a} \right)}
	+ \arccos \left( \sqrt{\frac{1}{a}} \right)
	+ \frac{\pi}{2}
	\right)
	\\
	&=\pi \sqrt{a^3} + 2\sqrt{a(a-1)}  + 2\sqrt{a^3} \arccos \left( \sqrt{\frac{1}{a}} \right),
	\\
	\sigma &= \sqrt{a} \arccos \left( \sqrt{\frac{1}{2} - \frac{a-1}{\sqrt{4a(a-1)}}} \right)
	\\
	&= \sqrt{a} \arccos \left( \sqrt{\frac{1}{2} - \frac{1}{2} \sqrt{\frac{a-1}{a}}} \right)
\end{align*}
and
\begin{align*}
	\mathcal{A}
	&=
	\left(
	\pi \sqrt{a^3} + 2\sqrt{a(a-1)}  + 2\sqrt{a^3} \arccos \left( \sqrt{\frac{1}{a}} \right)
	\right)
	\left(
	- \frac{1}{a} - \epsilon'
	\right)
	\\
	& \hphantom{=}
	+8 \sqrt{a} \arccos \left( \sqrt{\frac{1}{2} - \frac{1}{2} \sqrt{\frac{a-1}{a}}} \right).
\end{align*}
What happens here is that in fixed coordinates the angular velocity of~$M_3$ exactly matches that of~$M_2$ at the point of collision:
\begin{align} \label{eq angular velocity at pi/2}
	\frac{\diff \varphi}{\diff t} &=
	\frac{L}{r^2}
	=
	\sqrt{a\left(1-\epsilon^2\right)}
	=
	\sqrt{a\left(1-\frac{4a(a-1)}{4a^2}\right)}=1
\end{align}
Theorem~\ref{theorem Bolotin MacKay} can no longer exclude a collision in the perturbation to the restricted three-body problem because the ingoing and outgoing vectors in rotating coordinates are parallel.
So, at this point a new loop forms around~$M_2$.
More explicitly, in the restricted three-body problem, orbits coming from generating orbits with~$\theta<\pi/2$ will wind $I-1$ times around both~$M_1$ and~$M_2$ and then pass in between~$M_1$ and~$M_2$, while orbits coming from generating orbits with~$\theta>\pi/2$ will wind $I-1$ times around~$M_1$ and then wind once around~$M_2$ before starting over.
This transition from~$\theta=\pi$ and~$\epsilon'=-1$ to~$\theta=\pi$ and~$\epsilon'=+1$ via~$\theta=0$ happens in the~$\{ S_{-\alpha-1} \}$ segments in families~$f$ and ~$h$ as described in~\cite{Henon1}.
For generating orbits with~$I=1$ and corresponding continued orbits in the restricted three-body problem for small mass ratios as $\theta$ crosses~$\pi/2$ see also figure~\ref{figure crossing pi/2} in the next chapter.

\begin{figure}[t]
\centering
\begin{subfigure}[t]{0.45\textwidth}
	\centering
	\raisebox{2mm}{\includegraphics[width=\textwidth]{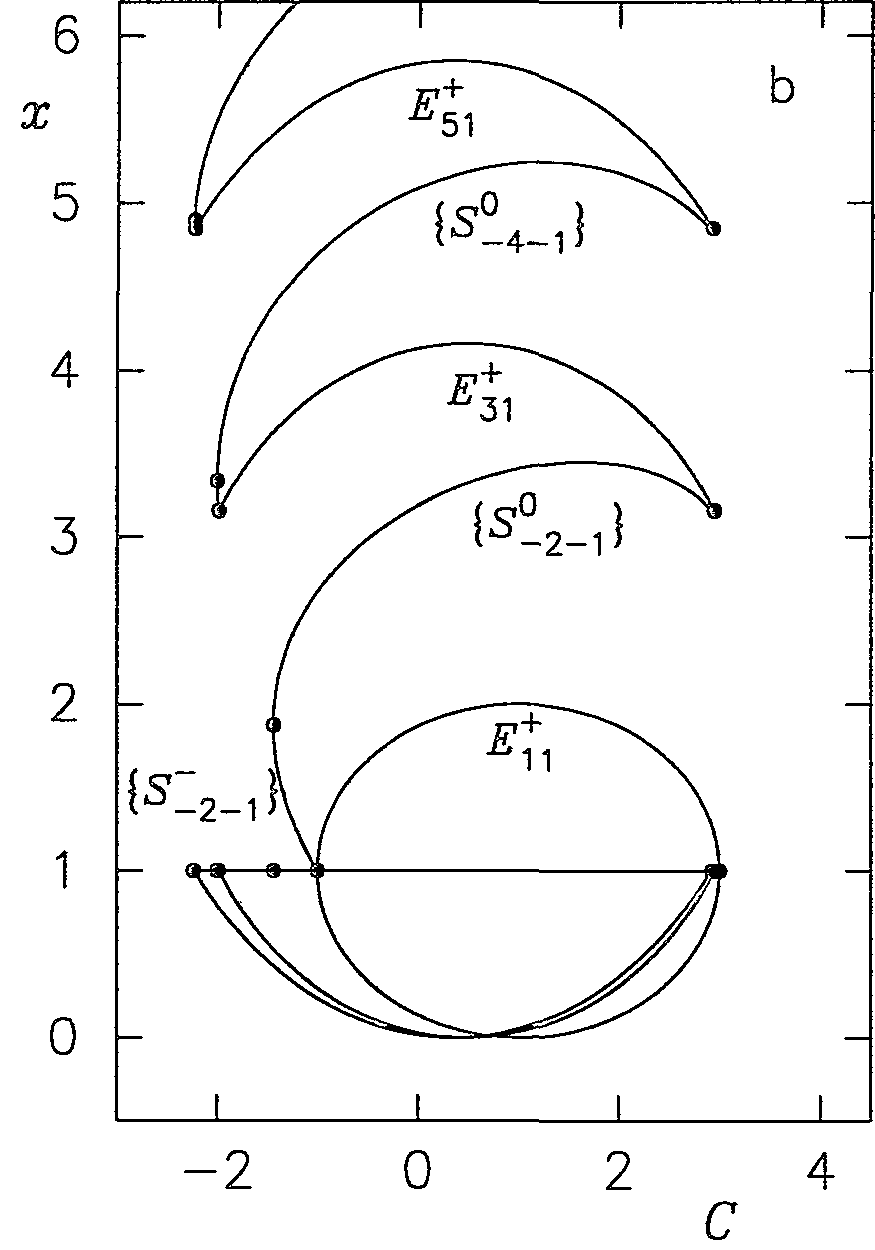}}
	\caption{Family~$f$}
	\label{figure family f}
\end{subfigure}
\hfill
\begin{subfigure}[t]{0.45\textwidth}
	\centering
	\includegraphics[width=\textwidth]{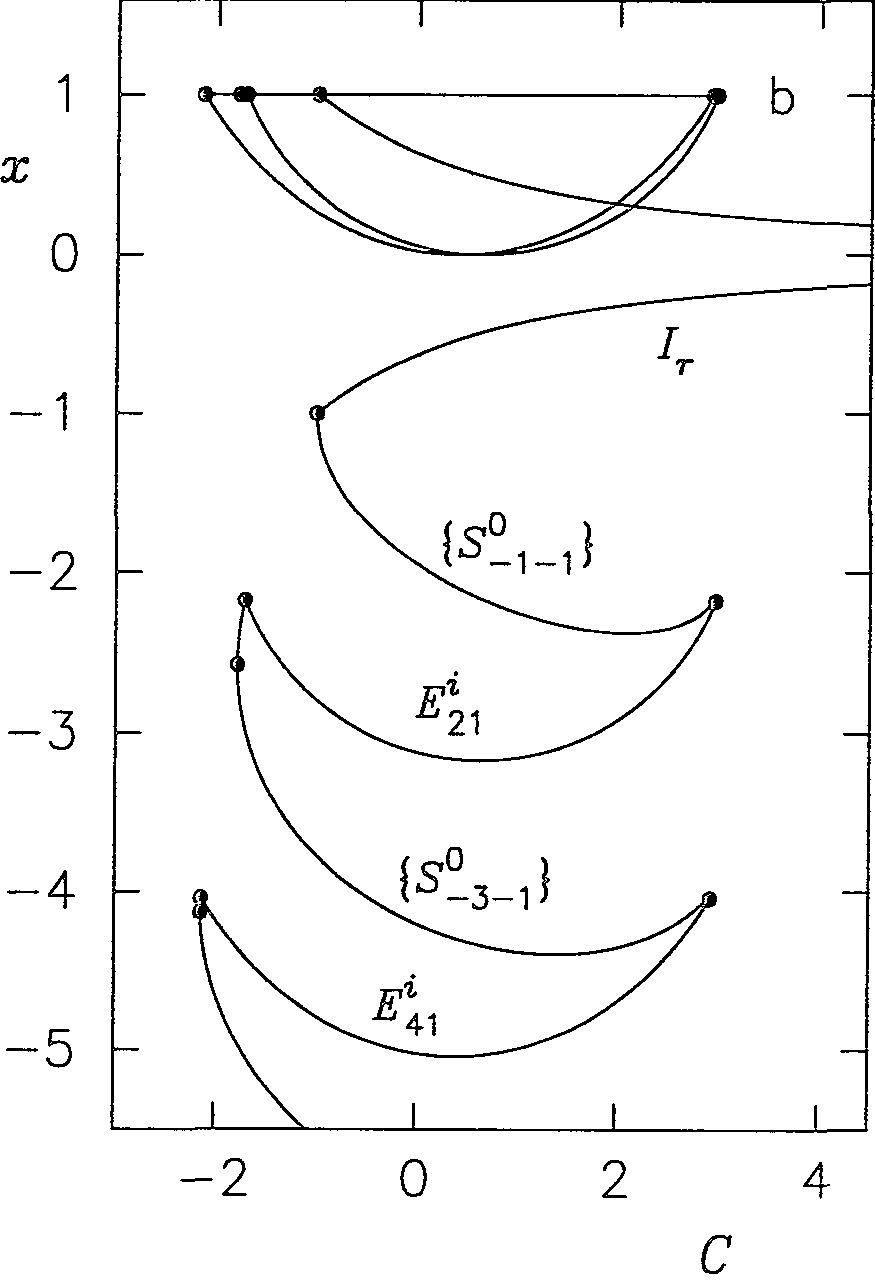}
	\caption{Family~$h$}
	\label{figure family h}
\end{subfigure}
\caption{Two symmetric families of generating orbits as shown in~\cite{Henon1}.}
\end{figure}

Generating family~$f$---see figure~\ref{figure family f}---comes from the simple direct orbits in Hill's lunar problem and transitions into family~$E_{1, 1}^+$ of first species generating orbits at the critical point with energy $H_0 = -3/2$ and abscissa~$q_1(0)=1$ for negative~$\dot{q}_2(0)$.
It follows~$E^+_{1, 1}$ with growing energy, undergoes collision with~$M_1$ at~$H_0 = - 1/2$ and~$q_1(0) = 2$, and picks up a loop around~$M_1$ there.
The generating family arrives at~$H_0 = 1/2$ and~$q_1(0)= 1$ as a double cover of the retrograde circular orbit, which is also a second species orbit with~$\theta = \pi$ and~$\epsilon'= -1$.
From there it follows the branch of second species generating orbits~$S_{-2, -1}$, transitioning to~$\theta=\pi$ and~$\epsilon'=+1$ and picking up a loop at~$\theta=\pi/2$ and~$\epsilon'=+1$ as described above.
The last orbit at~$\theta=\pi$ is again a first species orbit in the family~$E^+_{3 1}$ and this cycle continues indefinitely between~$E^+_{2k-1,1}$ and~$S_{-2k,-1}$ for all~$k>0$.

Family~$h$ comes from the retrograde circular orbits~$I_\mathrm{r}$ and then at~$a=1$ goes into the branch $S_{-1,-1}$ at~$\theta = \pi/2$ and~$\epsilon'=-1$.
Then it goes along, while taking up a loop, to~$\theta=\pi$ and~$\epsilon'=+1$ as described above.
There, it takes the branch of first species family~$E_{2 1}$, which takes up another loop at collision with~$M_1$, and carries on until the inner loop collides with~$M_2$ again.
From there on it resumes a similar pattern as family~$f$ and alternates between~$S_{-2k-1, -1}$ and~$E_{2k, 1}$ for all~$k>0$.

The summary of the information from this chapter which is important for the proof of the main theorem is
\begin{lemma}\label{lemma sequences negative energy}
	For every~$c \in [- \sqrt{2}, 0]$ there exists a sequence of second species generating orbits consisting of S-arcs with action tending towards negative infinity and energy converging strictly from below to~$c$.
\end{lemma}

%% file: main_proof.tex
We are now ready to put everything together and show up the consequences of these orbits for the existence of contact structures in the restricted three-body problem.
For this we first of all want to improve lemma~\ref{lemma sequences negative energy}.
In fact, we can get sequences of orbits with negative action tending towards negative infinity with constant energy at any value between~$-\sqrt{2}$ and 0.
We can also choose these sequences such that all orbits are ordinary and nondegenerate.
This is particularly helpful since the statement of theorem~\ref{theorem Bolotin MacKay} finds orbits close to the generating orbit with the exact same energy under these circumstances.
\begin{lemma} \label{lemma sequence with exact energy}
	For every~$c \in [-\sqrt{2}, 0)$ there exists a sequence of ordinary nondegenerate generating orbits of the second species with energy~$c$ and negative action tending to $- \infty$.
\end{lemma}
\begin{figure}[tb]
\centering
\includegraphics{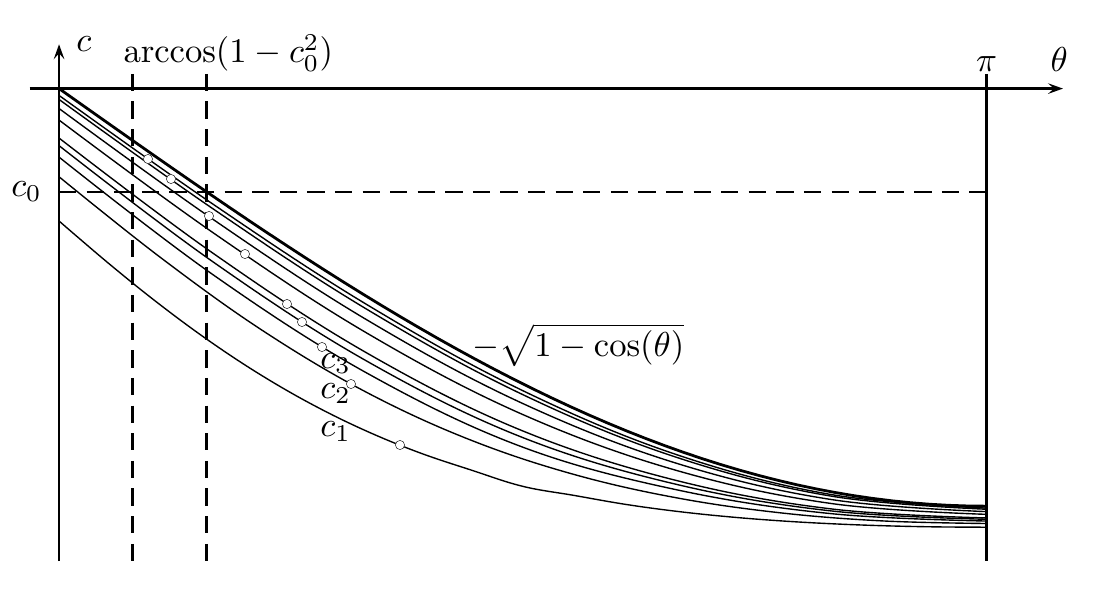}
\caption{The curves~$c_I$ mapping a simple direct arc of polar intersection angle~$\theta$ to its energy.}
\end{figure}

\begin{proof}
	As we have seen in the preceding chapter~\ref{chapter action of generating orbits}, there exists for every polar intersection angle~$\theta \in [0,\pi]$ with the unit circle and every number of rotation~$I \geq 1$ a unique second species arc.
	For~$\epsilon' = + 1$ these arcs have energy converging strictly monotonically to~$-\sqrt{1- \cos \theta}$ from below and their action tends to~$-\infty$ if~$\theta>0$ as~$I \to \infty$.
	Since the elapsed time~$\tau$ and the energy~$H_0 = c$ of the arcs depend smoothly on~$\theta$, we get for every~$I \geq 1$ a smooth curve
	\begin{align*}
		c_I \colon [0,\pi] &\to \RR
	\end{align*}
	mapping~$\theta$ to the energy of the unique direct arc where the timing condition~\eqref{eq timing condition} is satisfied for~$I$ and~$\theta$.
	The sequences for all fixed~$\theta$ are strictly increasing in~$c$ along~$a$, so we have~$c_{I_1} < c_{I_2}$ for all $I_1<I_2$ and the curves~$c_I$ converge to $-\sqrt{1-\cos \theta}$ as $I \to \infty$.
	
	Let~$c_0 \in [- \sqrt{2}, 0)$.
	We restrict the curves~$c_I$ to~$[\arccos(1-c_0^2)/2, \pi] \subset (0, \pi]$.
	Here, for all sequences of arcs with fixed~$\theta$, the action tends towards~$-\infty$ as~$a \to \infty$.
	So we can find for every~$\theta \in [\arccos(1-c_0^2)/2, \pi]$ an~$N_\theta>0$ such that the action of the arc with polar intersection angle~$\theta$ and semi-major axis~$a > N_\theta$ is negative.
	Assign this real number~$N_\theta$ to every~$\theta$ by the smooth function
	\begin{align*}
		N \colon [\arccos(1-c_0^2)/2, \pi] &\to \RR
		\\
		\theta & \mapsto N_\theta.
	\end{align*}
	Then~$N$ attains its maximum~$\hat{N}$ at some point and for all~$\theta \in [\arccos(1-c_0^2)/2, \pi]$ and~$a > \hat{N}$ the corresponding arcs have negative action.
	
	To find the smallest~$I$ where every point on the curve~$c_I$ is attained for an arc with semi-major axis~$a > \hat{N}$, we look at solutions of the timing condition for all real~$I_\theta \in (1, \infty)$.
	Again we assign this value to every arc with semi-major axis~$a = \hat{N}$ and~$\theta \in [\arccos(1-c_0^2)/2, \pi]$ by the smooth function
	\begin{align*}
		I \colon [\arccos(1-c_0^2)/2, \pi] &\to (1, \infty) \subset \RR
		\\
		\theta & \mapsto I_\theta
		.
	\end{align*}
	This, again, attains its maximum~$\hat{I} \in (0, \infty)$ at some point.
	So for all integers $I > \hat{I}$ we have full curves~$c_I \colon [\arccos(1-c_0^2)/2, \pi] \to (- \infty, 0)$ corresponding to arcs with semi-major axis~$a > \hat{N}$, \ie with negative action.
	Finally,
	\begin{align*}
		\left\{ \left( \theta, c_I(\theta) \right) \mid I > \hat{I}, \theta \in [\arccos(1-c_0^2)/2, \pi) \right\}
		\cap
		\left\{(\theta, c_0) \mid \theta \in [\arccos(1-c_0^2)/2, \pi)  \right\}
	\end{align*}
	is a sequence of points on the graphs~$\bigcup_I \Gamma(c_I)$ converging to $(\arccos(1-c_0^2), c_0)$.
	There are no degenerate or non-ordinary orbits in these sequences anymore since~$\theta \neq 0, \pi$.
	Furthermore, the only orbits with parallel velocities at collision in these sequences are where the arcs intersect the unit circle with angle~$\theta = \pi/2$.
	We can simply exclude these, since this affects at most one single element in the sequence.
	Hence, we have found a sequence of ordinary nondegenerate second species generating orbits with energy~$c_0$ and negative action.
	
	By doing the same process again for an arbitrarily small action we can show that the action of orbits in this sequence indeed tends towards~$-\infty$.
\end{proof}
In a next step we will check the integral of the first de Rham generator~$\beta_0$---which was computed in~\ref{lemma de Rham generator}---along the continued orbits that we get from the sequences of generating orbits.
\begin{figure}[t!]
	\centering
	\begin{subfigure}[b]{0.495\textwidth}
		\centering
		\includegraphics[width=\textwidth]{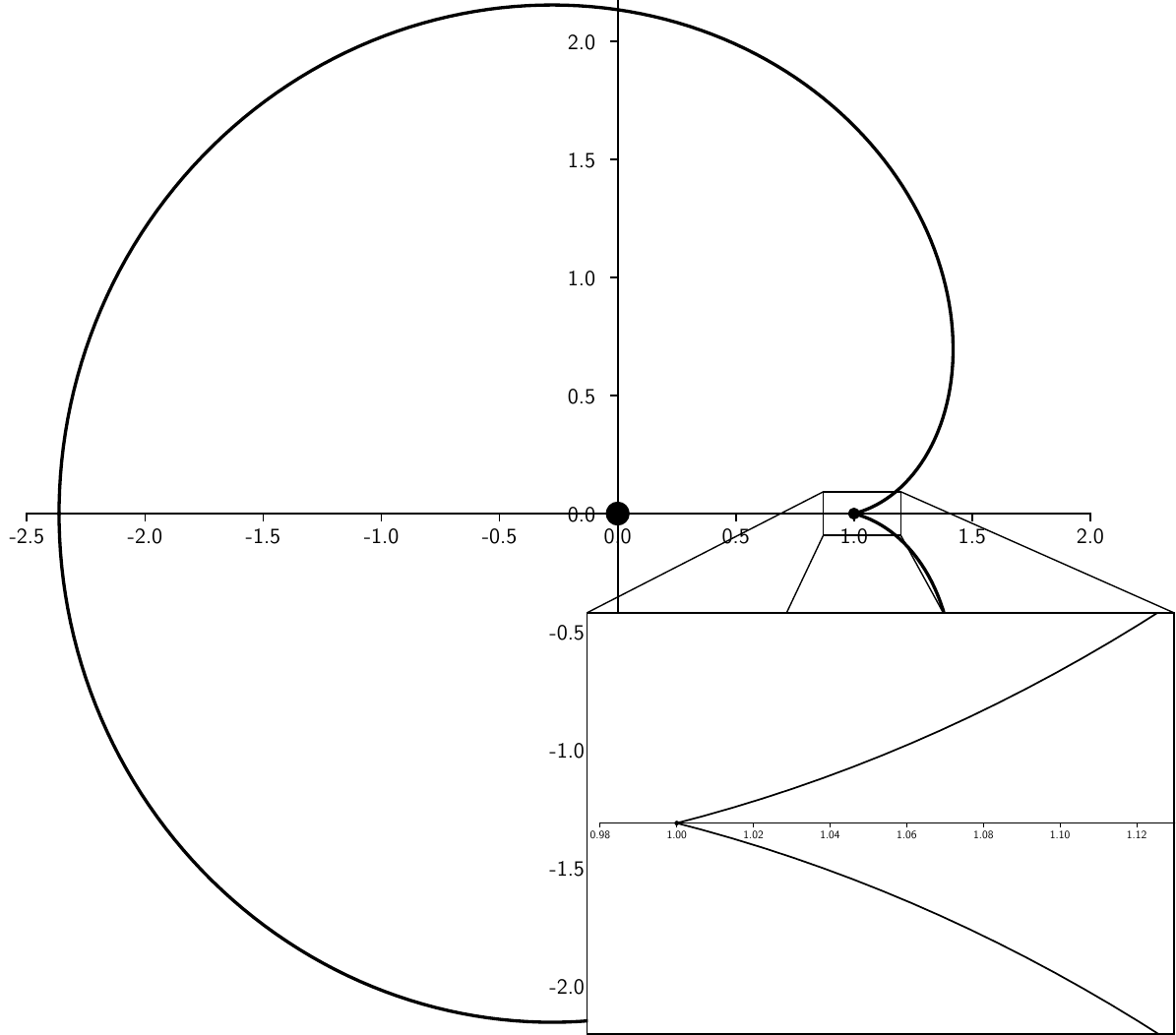}
		\caption{Generating orbit for $\theta = \pi/3$}
	\end{subfigure}
	\hfill
	\begin{subfigure}[b]{0.495\textwidth}  
		\centering 
		\includegraphics[width=\textwidth]{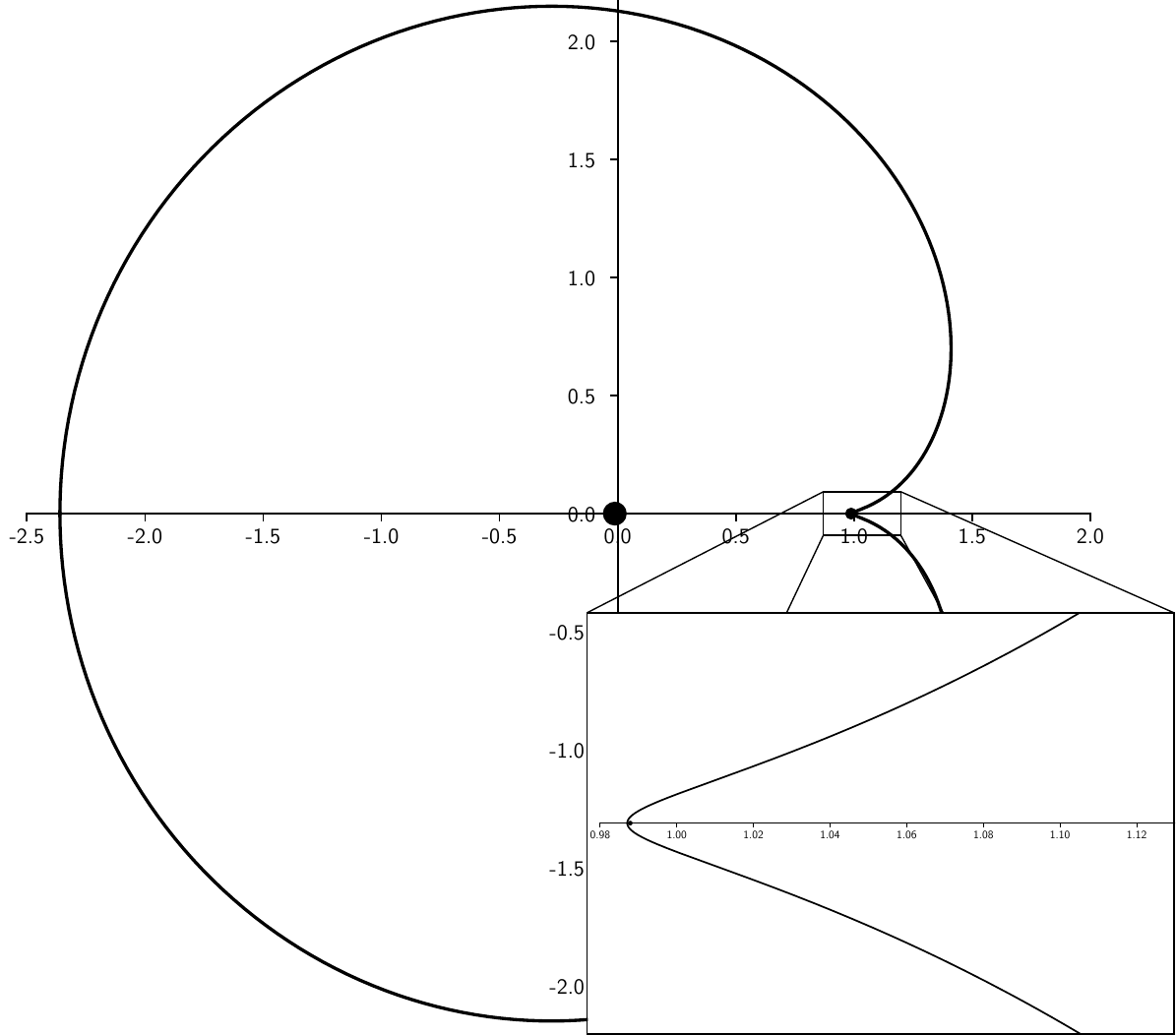}
		\caption{Continued orbit for $\theta = \pi/3$}
	\end{subfigure}
	\vskip\baselineskip
	\begin{subfigure}[b]{0.495\textwidth}   
		\centering 
		\includegraphics[width=\textwidth]{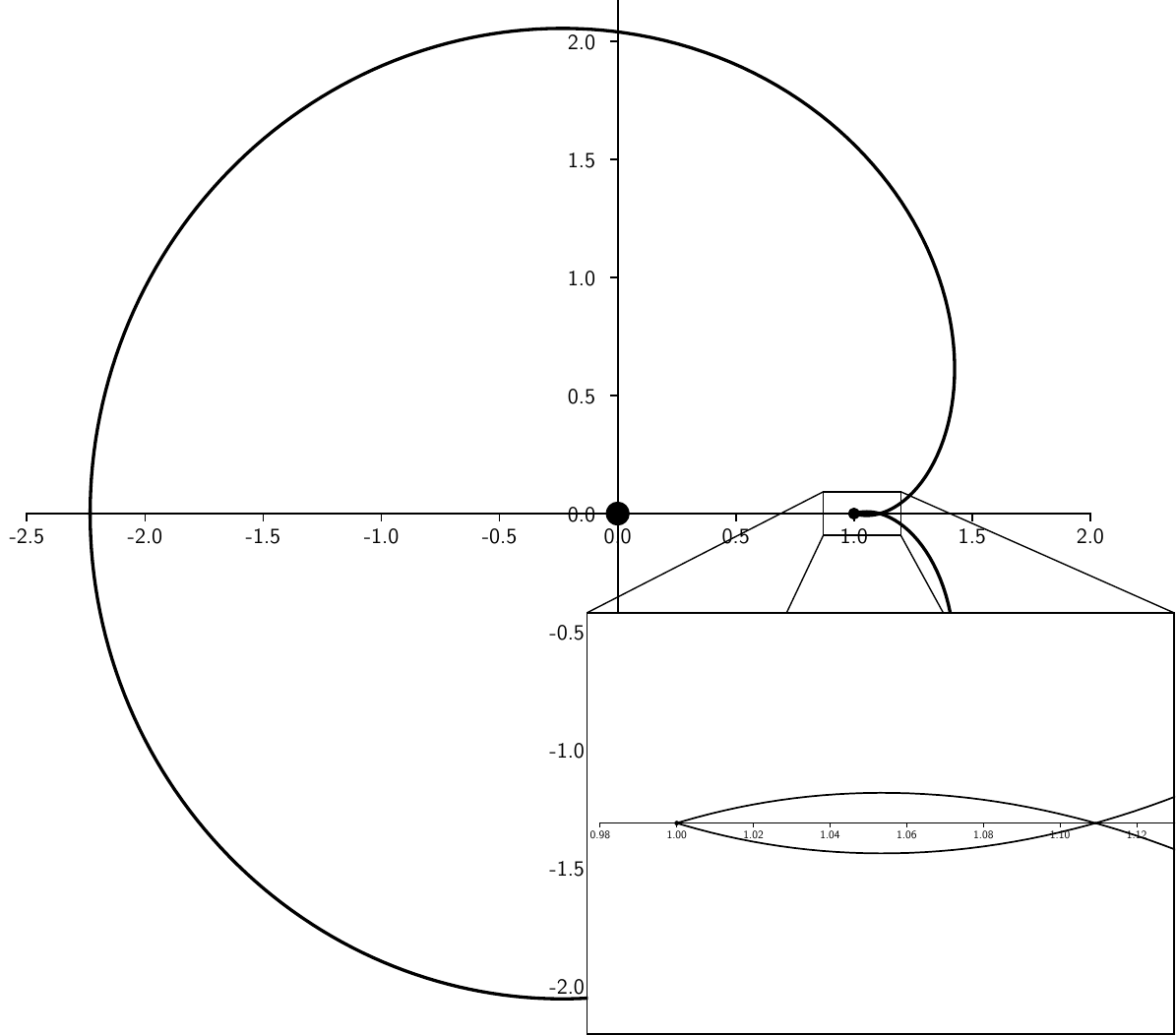}
		\caption{Generating orbit for $\theta = 2 \pi/3$}
	\end{subfigure}
	\hfill
	\begin{subfigure}[b]{0.495\textwidth}   
		\centering 
		\includegraphics[width=\textwidth]{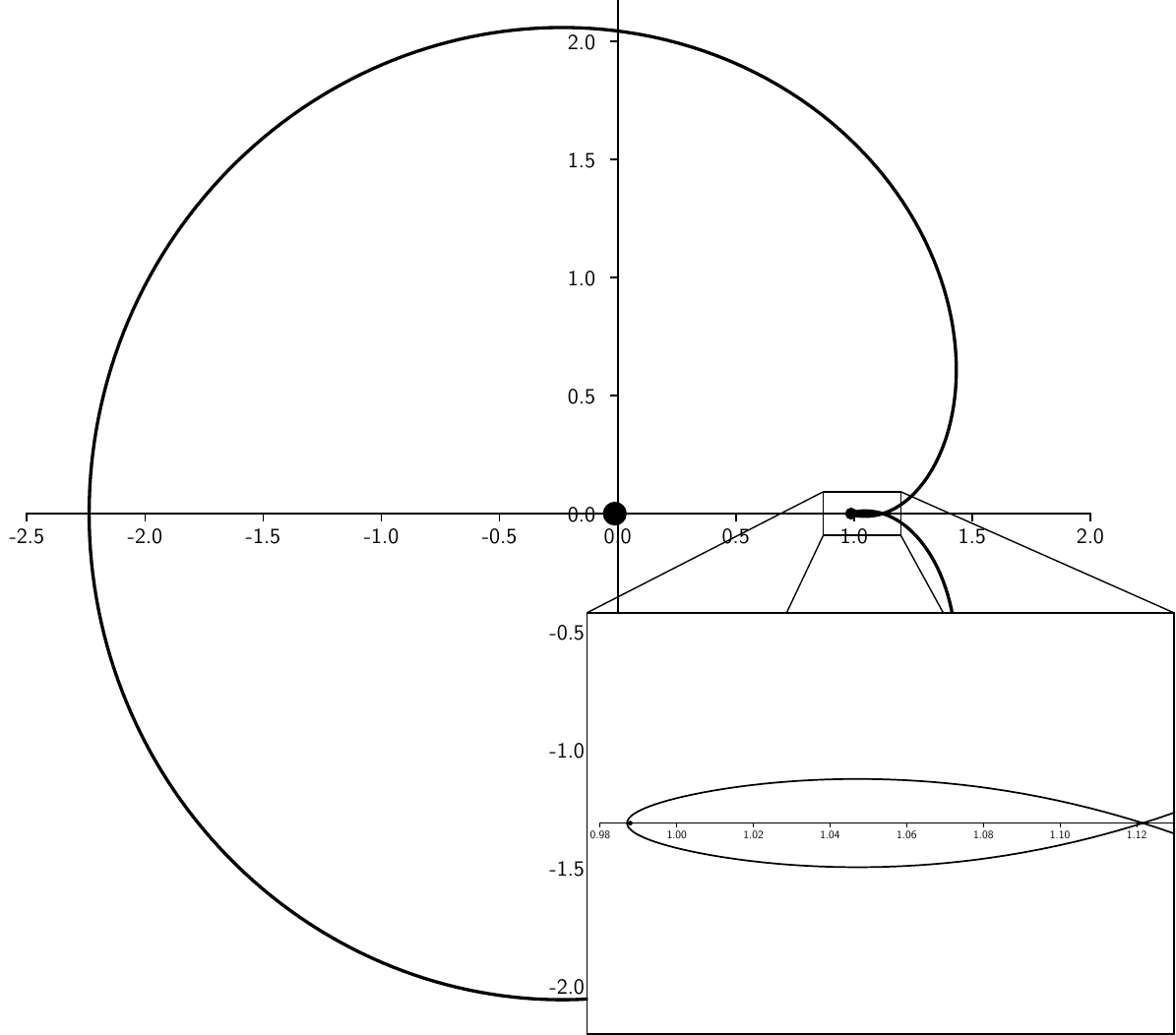}
		\caption{Continued orbit for $\theta = 2 \pi/3$}
	\end{subfigure}
	\caption{Generating orbits as~$\theta$ crosses~$\pi/2$ and their continuations at the same energy level in the restricted three-body problem with earth-moon mass ratio~$\mu = 0.01214$.}
	\label{figure crossing pi/2}
\end{figure}

\begin{lemma} \label{lemma integral of generator}
	Let~$\gamma$ be an orbit in one of the sequences from lemma~\ref{lemma sequence with exact energy} and~$\gamma_\mu$ a continuation to the restricted three-body problem with mass ratio~$\mu$ by theorem~\ref{theorem Bolotin MacKay}.
	Then the integral over the first de Rham generator~$\beta_0$ of~$\overline{\Sigma} \cong \overline{\Sigma}_{\mu, c}$ along~$\gamma_\mu$ is
	\begin{align*}
		\int_S^1 \gamma_\mu^*\beta_0 = -1.
	\end{align*}
\end{lemma}

\begin{proof}
	As in the proof of~\ref{lemma sequence with exact energy}, the initial and final velocities of the generating arc at collision are non-parallel in rotating coordinates and by theorem~\ref{theorem Bolotin MacKay} the continued orbit~$\gamma_\mu$ does not collide with~$M_2$.
	We can compute the integral over~$\beta_0 \vert_{\Sigma_{\mu, c}} = 2 \diff \vartheta - \diff \varphi_1 - \diff \varphi_2$ by two times the rotation number minus the two winding numbers around~$M_1$ and~$M_2$.
	
	In the case~$\theta \in (0, \pi/2)$ the angle between the initial and final velocities of the generating arc at collision is between 0 and~$\pi$ in rotating coordinates.
	The continuation uses hyperbolic solutions on the Levi-Civita regularisation close to~$M_2$ and all arcs are outgoing.
	So the continued orbit passes between~$M_1$ and~$M_2$ on a trajectory that looks like a hyperbola close to~$M_2$.
	Recall the integer numbers of rotation~$I$ and~$J$ of~$M_2$ and~$M_3$ around~$M_1$.
	In all of our sequences we have $J=0$ while $I \geq 1$.
	This gives us the rotation number~$-I$, the winding number~$-I$ around~$M_1$ and the winding number~$-I+1$ around~$M_2$.
	Overall, this gives
	\begin{align*}
		\int_{S^1} \gamma_\mu^*\beta_0 &= -2I + I + I - 1 = -1.
	\end{align*}
	In the case~$\theta \in (\pi/2, \pi)$ the angle between initial and final velocity at collision is in between~$\pi$ and~$2 \pi$, and the angular velocity dropping below one---see\eqref{eq angular velocity at pi/2}---forming an additional loop around~$M_2$ in the continuation.
	Therefore, the rotation number becomes~$-I-1$, the winding numbers~$-I$ around~$M_1$ and~$-I-1$ around~$M_2$.
	All in all we get the same integral
	\begin{align*}
		\int_{S^1} \gamma_\mu^* \beta_0 = 2(-I - 1) + I + I + 1 = -1,
	\end{align*}
	as claimed.
\end{proof}

Before we can finally prove the main theorem, we first need to discuss exactly what Hamiltonian structures and primitives we are dealing with:
\\
Let~$\overline{\Sigma}_{\mu, c}$ be the Moser-regularised energy hypersurface of the restricted three-body problem with mass ratio~$\mu \in (0,1)$ for an energy $c > H_\mu(L_5)$ above the highest critical value.
By the Moser regularisation, the Hamiltonian structure~$\omega$ on~$\Sigma_{\mu, c}$ extends to collisions and we get the Hamiltonian manifold~$(\overline{\Sigma}_{\mu, c}, \overline{\omega})$.
The Liouville 1-form~$\lambda = p \diff q$, on the other hand, does not extend.
We do, however, have a local primitive~$\tilde{\lambda}$ of~$\overline{\omega}$ in a neighbourhood~$U$ of collision at each primary because the fibre is Lagrangian and so~$\omega$ vanishes on the intersection of the regularised hypersurface with the fibre over~$p=\infty$.
On the intersection $V := \Sigma_{\mu, c} \cap U$ both~$\lambda \vert_V$ and~$\tilde{\lambda} \vert_V$ are primitives of~$\overline{\omega} \vert_V = \omega \vert_V$, so $\diff (\lambda  \vert_V - \tilde{\lambda}  \vert_V) = 0$.
The subset~$V \subset \overline{\Sigma}_{\mu, c}$ is diffeomorphic to~$S^*(\RR^n \setminus \{0\})$, which has trivial fundamental group for $n=3$.
So in this case the first de Rham cohomology of~$V$ is trivial and the closed 1-form~$\lambda \vert_V - \tilde{\lambda} \vert_V$ has a primitive~$f \in C^\infty (V)$.
Let~$g \in C^\infty (V)$ be a smooth cut-off function which is identically zero on~$\Sigma_{\mu, c}$ and identically one in a small neighbourhood of~$p= \infty$ in~$U$.
Define
\begin{align*}
	\overline{\lambda} :=
	\begin{cases}
		 \lambda & \text{in } \Sigma_{\mu, c} \setminus V
		 \\
		 \lambda - \diff(gf) & \text{in } V
		 \\
		 \tilde{\lambda} & \text{in } U \setminus V
	\end{cases}
\end{align*}
as a smooth extension of~$\lambda \vert_{\Sigma_{\mu, c} \setminus V}$ onto $\overline{\Sigma}_{\mu, c}$.
For~$n=2$, $V$ does have a large fundamental group, but we can simply define $\overline{\lambda}$ as the restriction of the three-dimensional construction.
The minimal distance to collision for all second species orbits with non-parallel velocity vectors at collision gives us enough space for our small neighbourhood~$U$.
In that way the local change of primitive does not change the action of the orbit.

Using this notation, we can now explicitly state the main theorem which will be divided into the planar case~$n=2$ and the spatial case~$n=3$:

\begin{theorem} \label{thm main theorem final}
	In the planar case we have:
	\\
	For every $c \in [ -\sqrt{2},0 )$ and~$r_0 \in \RR$ there exists a~$\mu_0 > 0$ such that for every~$\mu \in (0,\mu_0]$ the Hamiltonian manifold~$(\overline{\Sigma}_{\mu, c}, \overline{\omega})$ of the planar circular restricted three-body problem with mass ratio~$\mu$ does not admit an $\overline{\omega}$-compatible contact form~$\alpha$ such that~$[\alpha- \overline{\lambda}] = [r \beta_0] \in H^1_\mathrm{dR}(\overline{\Sigma}_{\mu, c})$ for any coefficient~$r \geq r_0$.
	\\
	For the spatial case:
	\\
	For every $c \in [ -\sqrt{2},0 )$ there exists a~$\mu_0 > 0$ such that for every~$\mu \in (0,\mu_0]$ the Hamiltonian manifold~$(\Sigma_{\mu, c}, \omega)$ of the circular restricted three-body problem with mass ratio~$\mu$ does not admit an $\omega$-compatible contact form~$\alpha$. 
\end{theorem}

\begin{proof}
	We begin with the planar case:
	Let~$c_0 \in [-\sqrt{2}, 0)$ and~$r_0 \in \RR$ be arbitrary.
	Then lemma~\ref{lemma sequence with exact energy} from above gives us an ordinary non-degenerate generating orbit~$\gamma_0$ with energy~$c$ and action~$\mathcal{A}(\gamma_0)< r_0 - \varepsilon$ for any~$\varepsilon>0$.
	By theorem~\ref{theorem Bolotin MacKay} there exists a~$\mu_0>0$ such that for all $\mu \in (0, \mu_0]$ there an orbit orbit~$\gamma_\mu$ in the restricted three-body problem with mass ratio~$\mu$ that has energy~$c$ and action~$\mathcal{A}(\gamma_\mu) \in (\mathcal{A}(\gamma_0) - \varepsilon, \mathcal{A}(\gamma_\mu) + \varepsilon)$, \ie $\mathcal{A}(\gamma_\mu)<r_0$.
	We have furthermore seen in lemma~\ref{lemma integral of generator} that the continued orbits from lemma~\ref{lemma sequence with exact energy} in the restricted three-body problem all have $\int \gamma_\mu^* \beta_0 = -1$.
	
	Assume by contradiction that we have an $\overline{\omega}$-compatible contact form~$\alpha \in \Omega^1(\overline{\Sigma}_{\mu, c})$ such that $[\alpha - \overline{\lambda}] = [r \beta_0] \in H_\mathrm{dR}^1(\overline{\Sigma})$ for an~$r \geq r_0$.
	Then we can write $\alpha = \overline{\lambda} + r \beta_0 + \diff f$  and integrate
	\begin{align*}
		\int_{S^1} \gamma_\mu^*\alpha
		=
		\int_{S^1} \gamma_\mu^* (\overline{\lambda} + r \beta_0 + \diff f)
		=
		\int_{S^1} \gamma_\mu^* \lambda + r \int_{S^1}\gamma_\mu^* \beta_0
		=
		\mathcal{A}(\gamma_\mu) - r
		\leq
		\mathcal{A}(\gamma_\mu) - r_0
		<
		0.
	\end{align*}
	This contradicts the fact that the integral over a compatible contact form along a Hamiltonian solution needs to be positive.
	
	In the spatial case all constructed orbits are contractible since already~$\Sigma_{\mu, c}$ has trivial fundamental group.
	Therefore, we directly get the claimed statement.
\end{proof}

This concludes the main part of this work and in the remaining chapters we will fist look at some numerical computations on these orbits and then discuss what questions arise and what future research could be suggested on this topic.

%% file: numerical.tex
After the analytical results from the previous chapter we want to add some numerical computations to visualise the continuations for actual mass ratios.
There will be two sections which both focus on direct orbits from sequences of section~\ref{subsection fixed theta}.
In the first section we will fix the number of rotation~$I$ of~$M_2$ around~$M_1$ at 1 and vary~$\theta$, \ie we will look at the first orbit in these sequences.
For the second section we will fix a~$\theta$ and look at the first several orbits in this sequence.

The data for generating orbits was obtained by numerically solving the timing condition~\eqref{eq timing condition} for the semi-major axis~$a$ and then computing the energy~$H_0$ and the action~$\mathcal{A}$ from equations~\eqref{eq energy fixed theta} and~\eqref{eq action fixed theta}.
The remaining parameter~$q_0$ represents the apoapsis distance from the origin, which is the initial position of the orbit.
Continued orbits for small mass ratios were found by using the symmetry~\eqref{eq reflection} and shooting perpendicularly from the $q_1$-axis while searching for perpendicular intersections when returning.
The programming for finding these orbits was done in python using standard libraries for solving ODEs and integration.
The results are numerical evidence of how far the described generating orbits might be followed and do not represent numerical or computer assisted proofs.

\subsection{First orbit in sequence with varying intersection angle}

\begin{figure}[tb]
\centering
\begin{subfigure}[b]{0.495\textwidth}
	\centering
	\includegraphics[width=\textwidth]{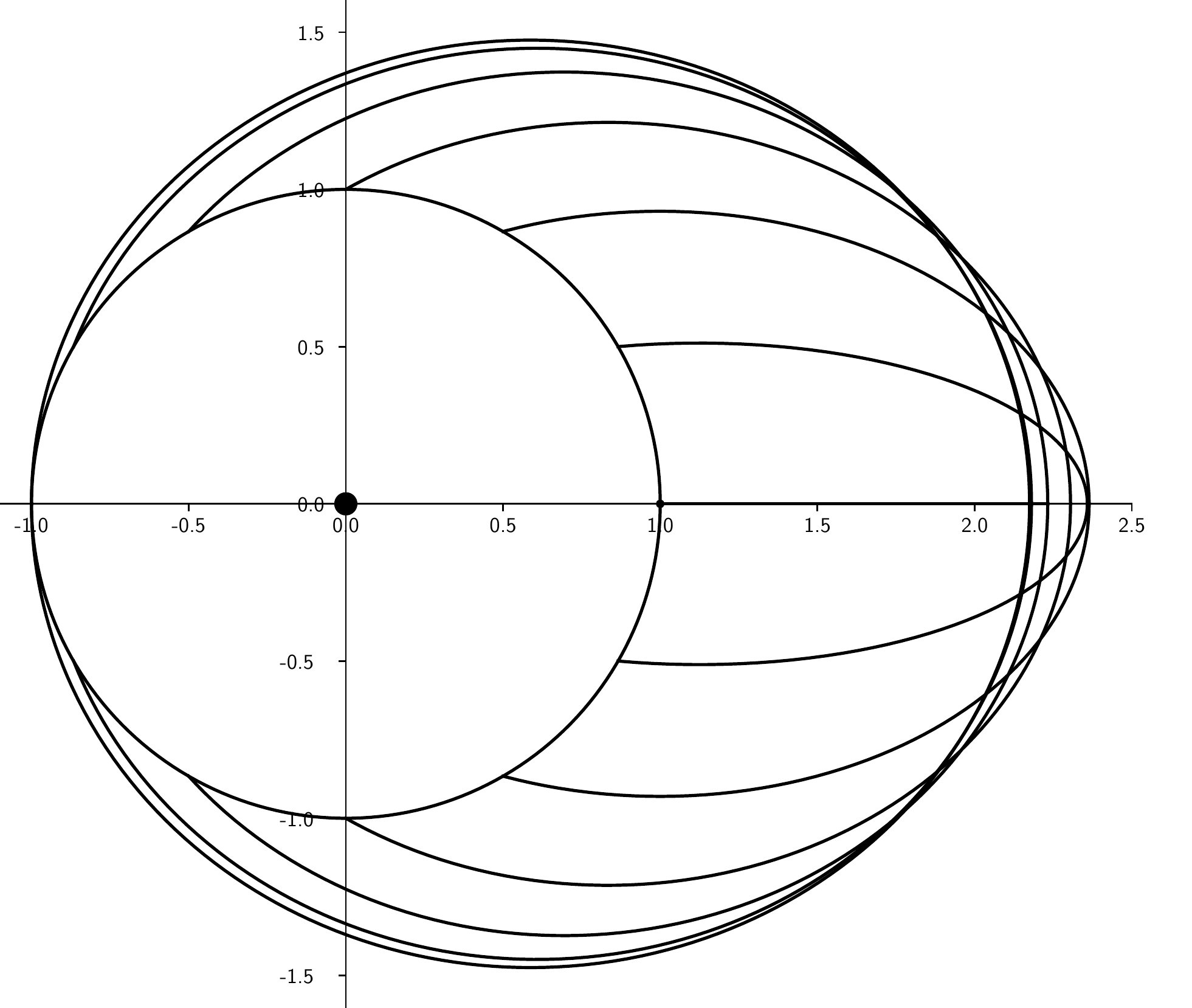}
	\caption{In fixed coordinates~$Q$.}
\end{subfigure}
\hfill
\begin{subfigure}[b]{0.495\textwidth}  
	\centering 
	\includegraphics[width=\textwidth]{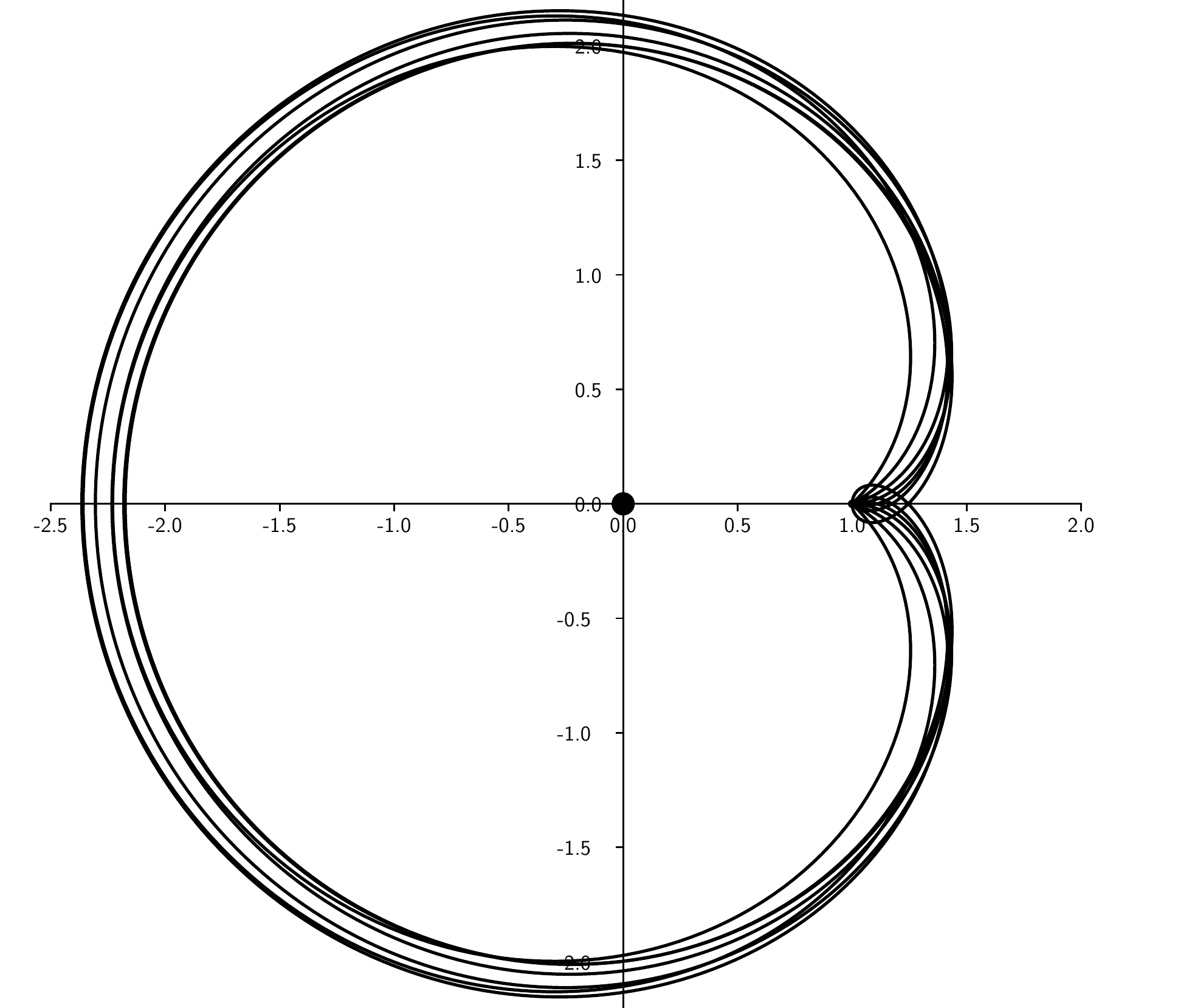}
	\caption{In rotating coordinates~$q$.}
\end{subfigure}
\vskip\baselineskip
\begin{subfigure}[b]{0.495\textwidth}   
	\centering 
	\includegraphics[width=\textwidth]{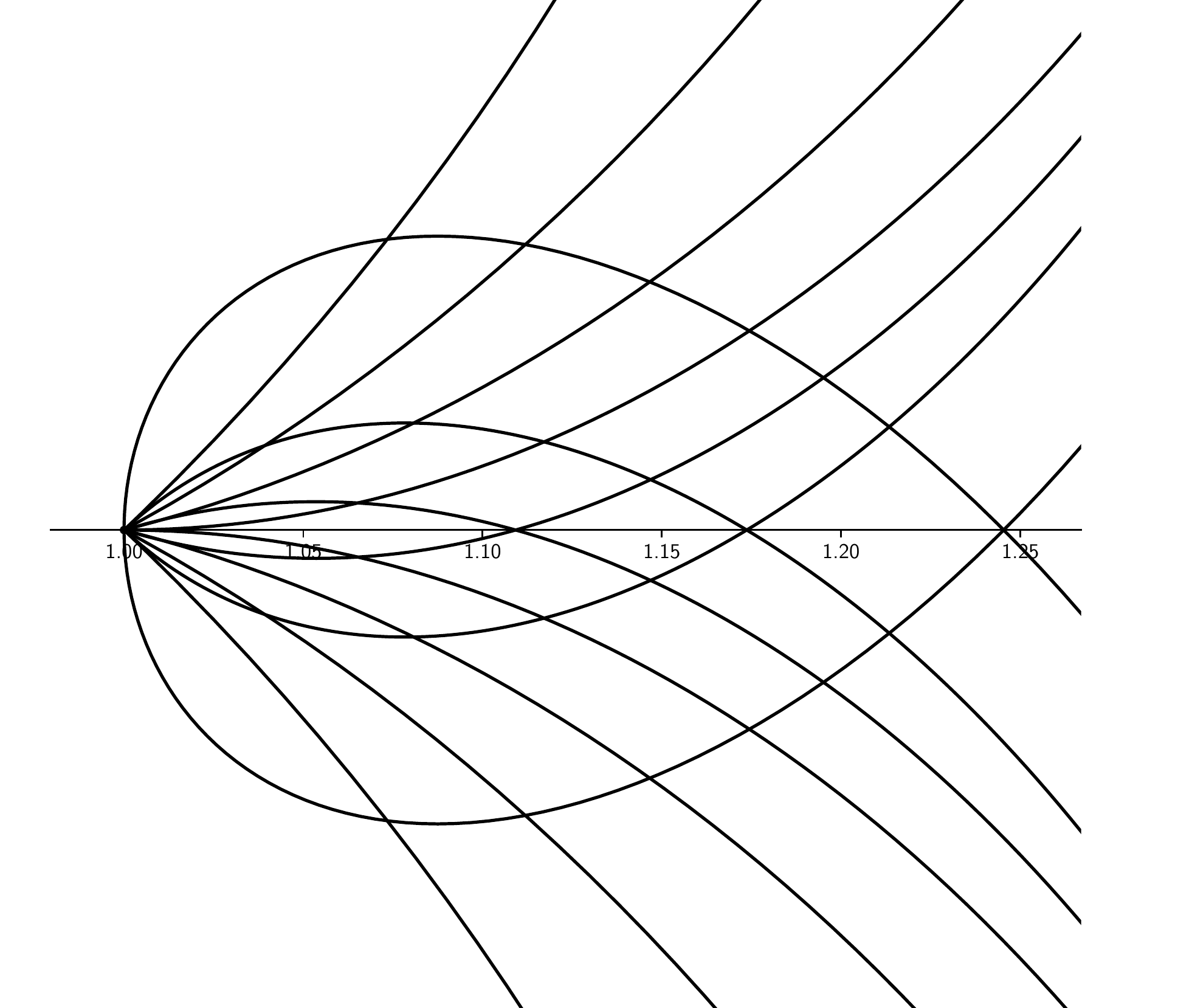}
	\caption{Close-up view of the neighbourhood of~$M_2$ in rotating coordinates.}
\end{subfigure}
\caption{Generating orbits for~$I=1$ and varying $\theta$ from 0 to $\pi$ in steps of 10 degrees.}
\label{figure generating all theta I1}
\end{figure}

Table~\ref{table generating all theta I1} shows the data for generating orbits for~$I=1$ and varying $\theta$ from 0 to ~$\pi$.
We use degrees instead of radians to describe the angle~$\theta$ which will be incremented in 10 degree steps.
From~\eqref{eq action fixed theta} we know that the action of generating orbits in the sequences with fixed~$\theta \in (0, \pi]$ tends to~$-\infty$ as~$I \to \infty$.
However, not all sequences always have negative action, as one can see in the case where $\theta$ is at 10 degrees.
In the boundary case~$\theta = 0$ the action of orbits in the sequence does not tend to negative values and is instead always positive.
Figure~\ref{figure generating all theta I1} shows every third of these generating orbits in fixed coordinates~$Q$ and in rotating coordinates~$q$.
We note the special cases~$\theta = \pi$ where the deflection angle is zero, \ie the generating orbit is non-ordinary, and $\theta = \pi/2$ where the deflection angle is~$-\pi$ and even the continued orbit might collide with~$M_2$.

Table~\ref{table continued astro all theta I1} shows the initial positions of these orbits when continued to the astronomical mass ratios $\mu \approx 0.000953$ of the Sun-Jupiter system, $\mu \approx 0.012143$ of the Earth-Moon system and $\mu \approx 0.108511$ of the Pluto-Charon system.
These mass ratios were chosen for the relevance in our solar system: The Sun-Jupiter system has the largest mass ratio of all the planets in relation to the sun and governs many phenomena like the occurrence of asteroids.
Obviously, the Earth-Moon system is of major importance to us and plays a large part in near-Earth satellites and lunar space missions.
Pluto and Charon have the largest mass ratio of a binary system in our solar system and serves as the largest mass ratio we will have to deal with in our current reach.
Table~\ref{table continued nonastro all theta I1} adds some larger non-astronomical values as well.
In all the continuations the same energy value is held as the generating orbit in accordance with theorem~\ref{theorem Bolotin MacKay}.

We can see that smaller values of~$\theta$ seem to continue up to higher mass ratios as opposed to larger values of~$\theta$.
Here, after some point no orbit could be found nearby and the corresponding entries are marked with an \anfzeichen{x}.
This makes good sense since the limit orbits at $\theta = \pi$ are bifurcation orbits and here two intersecting families of generating orbits split in the continuation and move away from the original bifurcation orbit.

Also the action of continued orbits seems to always be slightly larger than that of the generating orbit.
On the other hand, even the first generating orbits in most of theses sequences yield negative action, especially for large values of~$\theta$.
The smallest energy where one gets a generating orbit with negative action is already very close to the critical value~$-3/2$.
Although these are only isolated orbits, one might now expect that there is an obstruction for the restricted three-body problem to be of contact type for all energies between 0 and the highest critical value~$H_\mu(L_5)$.

\subsection{Sequence at 10 degrees intersection angle}

We now take a closer look at the first sequence from~\ref{table generating all theta I1} where the intersection angle~$\theta$ is nonzero.
This was the only one where the action of the first orbit is nonnegative.
So, we are firstly interested in how quickly the action becomes negative, but also if we can still continue the generating orbit to the large mass ratio of~$\mu = 0.5$ as the sequence goes on.

\begin{figure}[tb]
\centering
\begin{subfigure}[b]{0.495\textwidth}
	\centering
	\includegraphics[width=\textwidth]{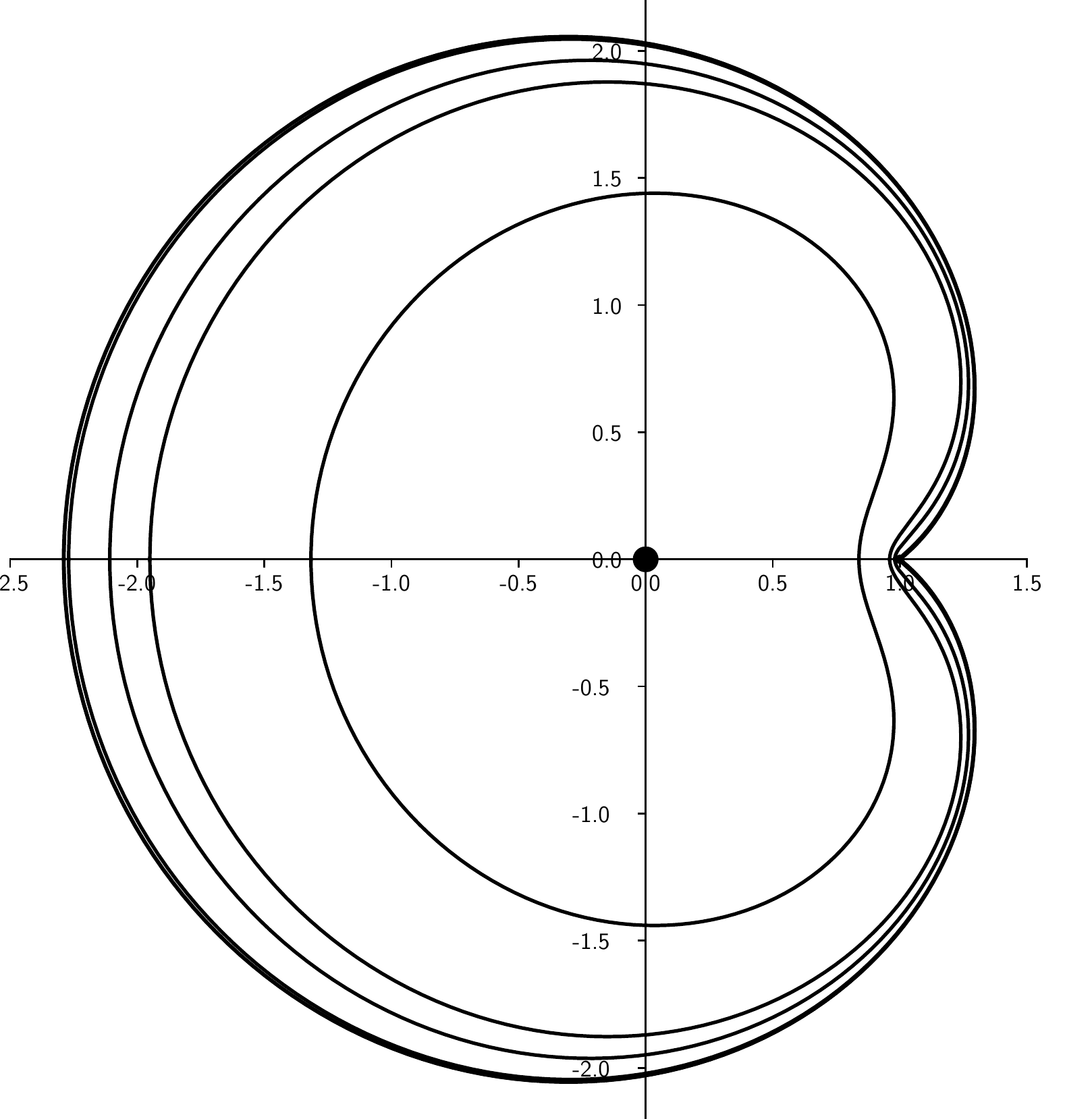}
	\caption{$I=1$}
\end{subfigure}
\hfill
\begin{subfigure}[b]{0.495\textwidth}  
	\centering 
	\includegraphics[width=\textwidth]{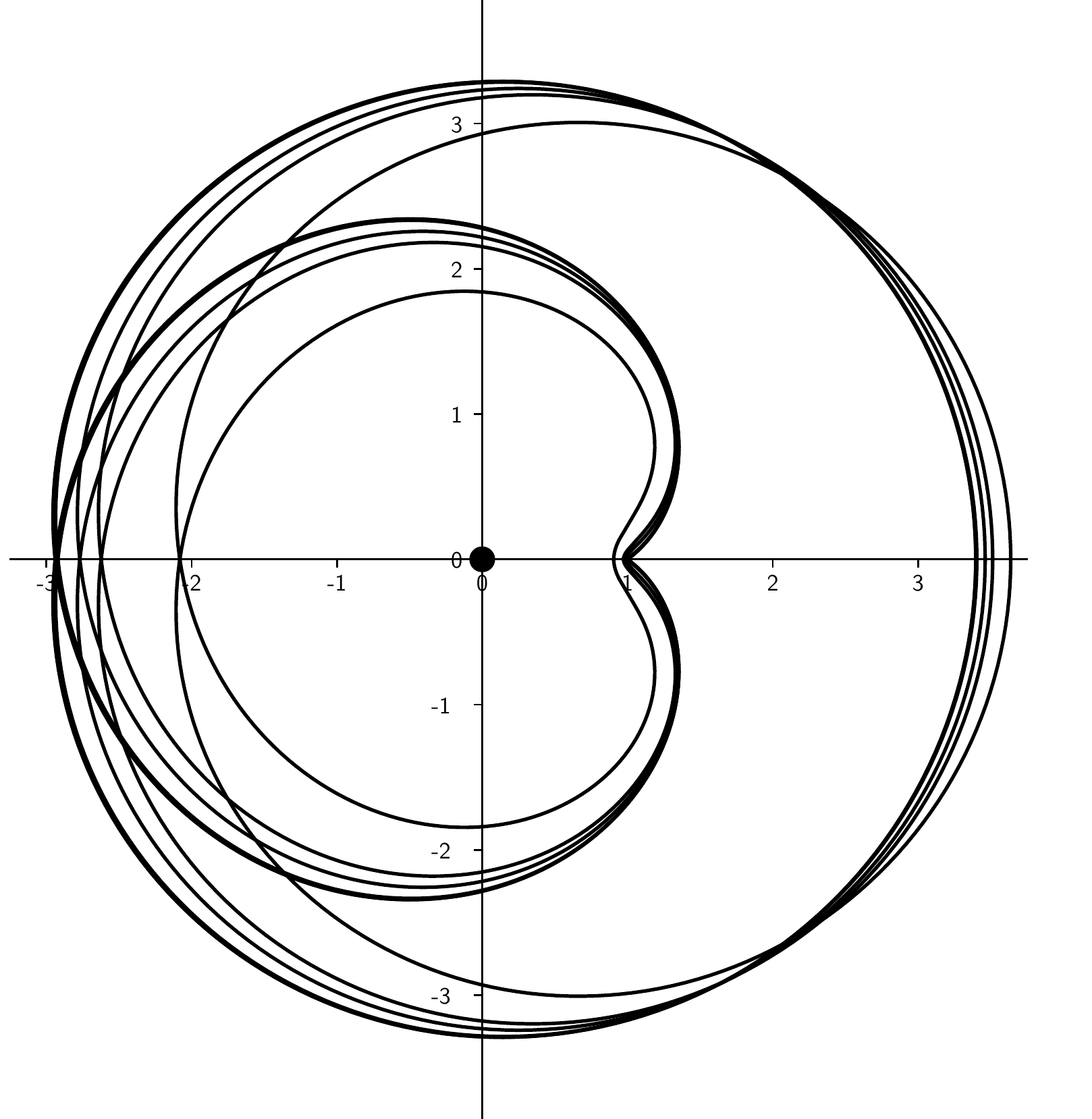}
	\caption{$I=2$}
\end{subfigure}
\vskip\baselineskip
\begin{subfigure}[b]{0.495\textwidth}   
	\centering 
	\includegraphics[width=\textwidth]{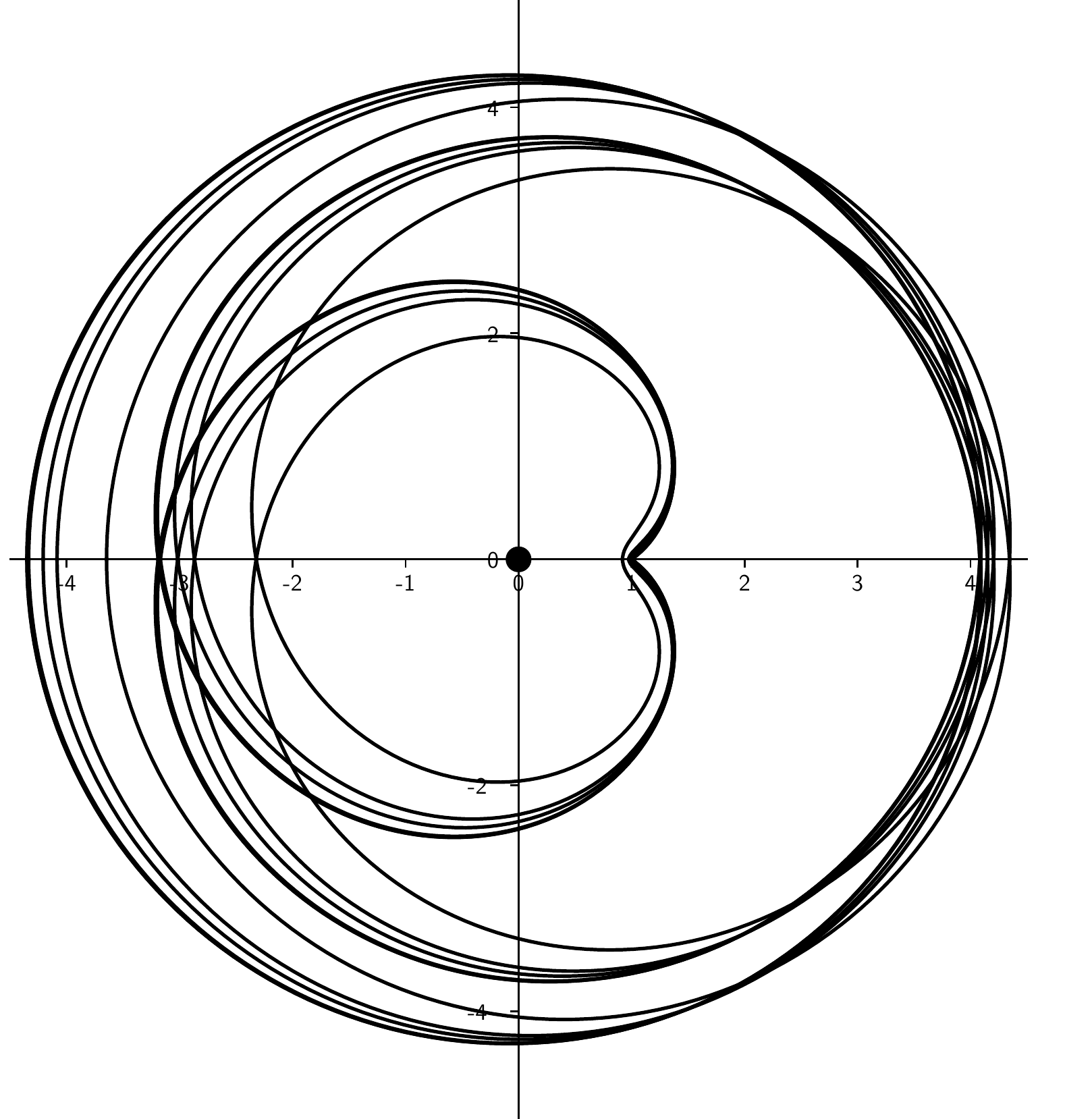}
	\caption{$I=3$}
\end{subfigure}
\hfill
\begin{subfigure}[b]{0.495\textwidth}   
	\centering 
	\includegraphics[width=\textwidth]{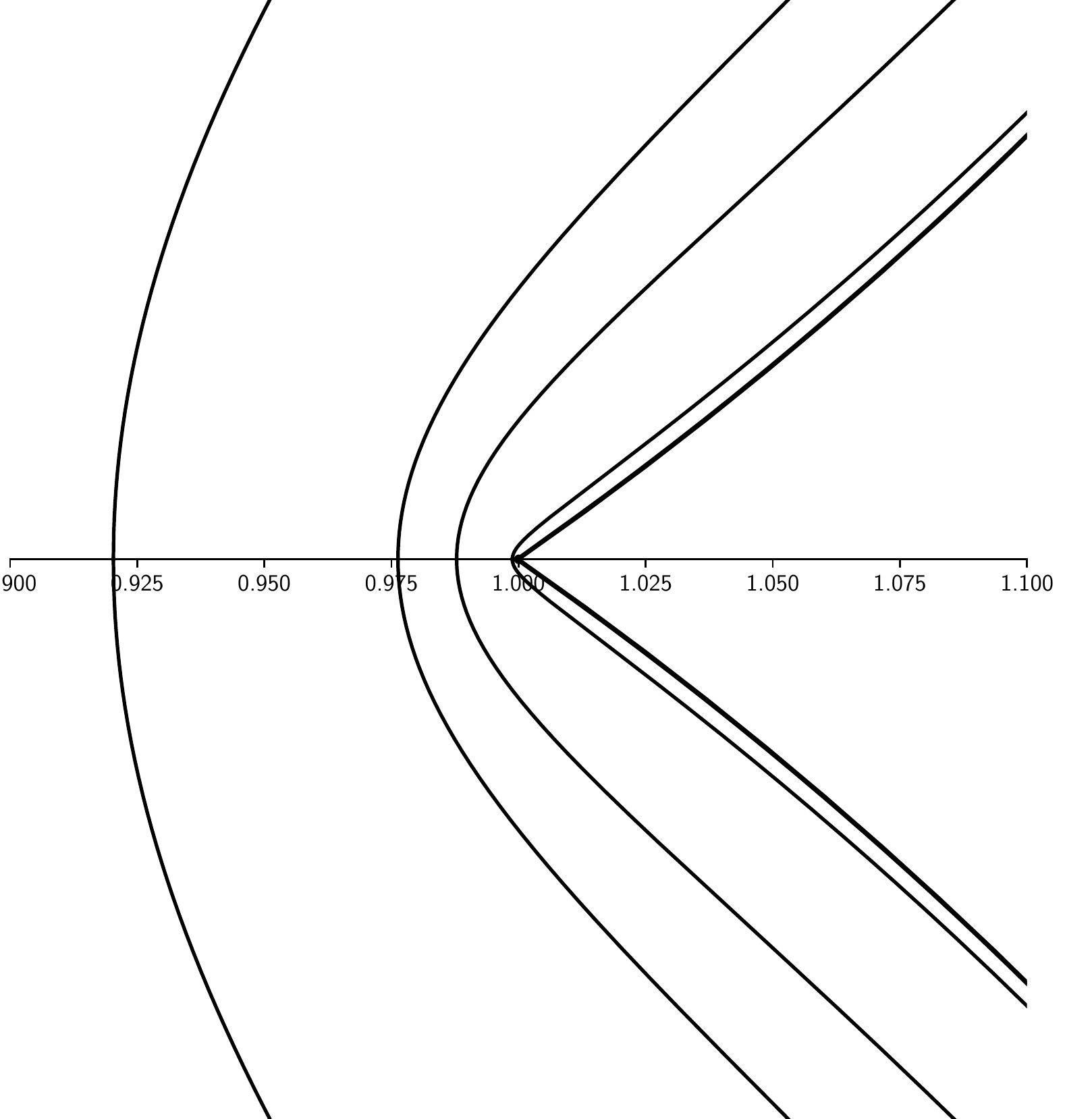}
	\caption{Close-up view of the neighbourhood of~$M_2$.}
\end{subfigure}
\caption{Generating orbit and continued orbits for $\theta = \pi/18$.}
\label{figure theta10}
\end{figure}

These generating orbits are computed in table~\ref{table generating theta10} and we see that the action becomes negative at the 10\textsuperscript{th} orbit in the sequence.
Furthermore, we can see the energy slowly converging to~$-\sqrt{1 - \cos (\pi/18)} \approx -0.123256$ while the semi-major axis grows.
Tables~\ref{table continued astro theta10} and~\ref{table continued nonastro theta10} show again the continued orbits to the same mass ratios as before.
Here, we notice that also the higher orbits in the sequence so far all continue all the way up to~$\mu = 0.5$ and eventually even the action becomes negative for this large mass ratio.
This evidence suggests that the restricted three-body problem might not be of contact type in between the highest critical value and zero for any mass ratio~$\mu \in (0,1)$.
Note, as mentioned in section~\ref{subsection fixed theta}, that all orbits with even~$I$ belong to the same family~$f$ and all orbits with odd~$I$ belong to family~$h$ in Hénon's notation.
This is visualised in figure~\ref{figure theta10} where the generating and continued orbits of this sequence are depicted for the first three~$I = 1,2,3$.
\clearpage

\begin{table}[tb]
	\centering
	\begin{tabular}{|c||c|c|c|c|} 
		\hline
		$\theta$ & $a$ & $q_0$ & $H_0$ & $\mathcal{A}$  \\ 
		\hline\hline
		0                               & 1.114891                 & 2.229783                   & -0.448474                  & 1.434641                            \\ 
		\hline
		10                              & 1.151460                 & 2.289513                   & -0.597508                  & 0.472433                            \\ 
		\hline
		20                              & 1.191803                 & 2.331984                   & -0.737353                  & -0.479144                           \\ 
		\hline
		30                              & 1.234738                 & 2.358649                   & -0.865060                  & -1.392612                           \\ 
		\hline
		40                              & 1.278895                 & 2.371414                   & -0.978834                  & -2.246054                           \\ 
		\hline
		50                              & 1.322831                 & 2.372473                   & -1.077948                  & -3.024066                           \\ 
		\hline
		60                              & 1.365152                 & 2.364145                   & -1.162569                  & -3.717803                           \\ 
		\hline
		70                              & 1.404646                 & 2.348721                   & -1.233529                  & -4.324279                           \\ 
		\hline
		80                              & 1.440377                 & 2.328356                   & -1.292090                  & -4.845201                           \\ 
		\hline
		90                              & 1.471746                 & 2.304988                   & -1.339732                  & -5.285607                           \\ 
		\hline
		100                             & 1.498494                 & 2.280307                   & -1.377985                  & -5.652556                           \\ 
		\hline
		110                             & 1.520668                 & 2.255753                   & -1.408307                  & -5.954007                           \\ 
		\hline
		120                             & 1.538560                 & 2.232546                   & -1.432015                  & -6.197965                           \\ 
		\hline
		130                             & 1.552639                 & 2.211726                   & -1.450243                  & -6.391878                           \\ 
		\hline
		140                             & 1.563480                 & 2.194211                   & -1.463929                  & -6.542246                           \\ 
		\hline
		150                             & 1.571719                 & 2.180859                   & -1.473821                  & -6.654366                           \\ 
		\hline
		160                             & 1.578013                 & 2.172531                   & -1.480481                  & -6.732160                           \\ 
		\hline
		170                             & 1.583022                 & 2.170155                   & -1.484294                  & -6.777997                           \\ 
		\hline
		180                             & 1.587401                 & 2.174802                   & -1.485467                  & -6.792454                           \\
		\hline
	\end{tabular}
	\caption{Generating orbits for $I = 1$ and ~$\theta$ in degrees.}
	\label{table generating all theta I1}
\end{table}

\begin{table}[tb]
	\centering
	\begin{tabular}{|c||c|c|c|c|c|c|} 
		\hline
		& \multicolumn{2}{c|}{Sun-Jupiter} & \multicolumn{2}{c|}{Earth-Moon} & \multicolumn{2}{c|}{Pluto-Charon}  \\ 
		\hline
		$\theta$ & $q_0$     & $\mathcal{A}$                & $q_0$     & $\mathcal{A}$               & $q_0$     & $\mathcal{A}$                  \\ 
		\hline\hline
		0     & 2.229061 & 1.446733              & 2.220504 & 1.546640             & 2.142401 & 2.136580                \\ 
		\hline
		10    & 2.288901 & 0.485253              & 2.281639 & 0.590492             & 2.216384 & 1.209834                \\ 
		\hline
		20    & 2.331460 & -0.465552             & 2.325234 & -0.355078            & 2.269946 & 0.288979                \\ 
		\hline
		30    & 2.358198 & -1.378186             & 2.352824 & -1.262459            & 2.305515 & -0.597382               \\ 
		\hline
		40    & 2.371026 & -2.230715             & 2.366375 & -2.109609            & 2.325631 & -1.426317               \\ 
		\hline
		50    & 2.372143 & -3.007720             & 2.368130 & -2.881020            & 2.332950 & -2.181565               \\ 
		\hline
		60    & 2.363871 & -3.700338             & 2.360452 & -3.567767            & 2.330171 & -2.853621               \\ 
		\hline
		70    & 2.348508 & -4.305571             & 2.345684 & -4.166823            & 2.319969 & -3.439064               \\ 
		\hline
		80    & 2.328214 & -4.825112             & 2.326039 & -4.679890            & 2.305002 & -3.939387               \\ 
		\hline
		90    & 2.304939 & -5.263994             & 2.303542 & -5.112050            & 2.288258 & -4.359638               \\ 
		\hline
		100   & 2.280390 & -5.629270             & 2.280016 & -5.470452            & 2.281012 & -4.707251               \\ 
		\hline
		110   & 2.256040 & -5.928906             & 2.257156 & -5.763208            & x        & x                       \\ 
		\hline
		120   & 2.233168 & -6.170932             & 2.236858 & -5.998577            & x        & x                       \\ 
		\hline
		130   & 2.212958 & -6.362855             & x        & x                    & x        & x                       \\ 
		\hline
		140   & 2.196764 & -6.511324             & x        & x                    & x        & x                       \\ 
		\hline
		150   & 2.190318 & -6.622143             & x        & x                    & x        & x                       \\ 
		\hline
		160   & x        & x                     & x        & x                    & x        & x                       \\ 
		\hline
		170   & x        & x                     & x        & x                    & x        & x                       \\ 
		\hline
		180   & x        & x                     & x        & x                    & x        & x                       \\
		\hline
	\end{tabular}
\caption{Continued orbits for $I=1$ and~$\theta$ in degrees for astronomical mass ratios.}
\label{table continued astro all theta I1}
\end{table}

\begin{table}[tb]
	\centering
	\begin{tabular}{|c||c|c|c|c|} 
		\hline
		& \multicolumn{2}{c|}{$\mu=0.2$} & \multicolumn{2}{c|}{$\mu=0.5$}  \\ 
		\hline
		$\theta$ & $q_0$     & $\mathcal{A}$           & $q_0$     & $\mathcal{A}$            \\ 
		\hline\hline
		0     & 2.059605 & 2.593750         & x        & x                 \\ 
		\hline
		10    & 2.150085 & 1.693790         & 1.816536 & 3.069224          \\ 
		\hline
		20    & 2.215887 & 0.792632         & 1.992154 & 2.275666          \\ 
		\hline
		30    & 2.261108 & -0.079578        & 2.113013 & 1.462972          \\ 
		\hline
		40    & 2.289277 & -0.898644        & x        & x                 \\ 
		\hline
		50    & 2.303809 & -1.647264        & x        & x                 \\ 
		\hline
		60    & 2.308293 & -2.315120        & x        & x                 \\ 
		\hline
		70    & 2.307463 & -2.898274        & x        & x                 \\ 
		\hline
		80    & x        & x                & x        & x                 \\ 
		\hline
		90    & x        & x                & x        & x                 \\ 
		\hline
		100   & x        & x                & x        & x                 \\ 
		\hline
		110   & x        & x                & x        & x                 \\ 
		\hline
		120   & x        & x                & x        & x                 \\ 
		\hline
		130   & x        & x                & x        & x                 \\ 
		\hline
		140   & x        & x                & x        & x                 \\ 
		\hline
		150   & x        & x                & x        & x                 \\ 
		\hline
		160   & x        & x                & x        & x                 \\ 
		\hline
		170   & x        & x                & x        & x                 \\ 
		\hline
		180   & x        & x                & x        & x                 \\
		\hline
	\end{tabular}
\caption{Contiunued orbits for $I=1$ and $\theta$ in degrees for non-astronomical mass ratios.}
\label{table continued nonastro all theta I1}
\end{table}

\begin{table}[tb]
	\centering
	\begin{tabular}{|c||c|c|c|c|} 
		\hline
		$I$  & $a$        & $q_0$      & $H_0$    & $\mathcal{A}$\\ 
		\hline\hline
		1  & 1.151461 & 2.289514  & -0.597508 & 0.472433   \\ 
		\hline
		2  & 1.703950 & 3.397169  & -0.439704 & 0.971722   \\ 
		\hline
		3  & 2.180544 & 4.351247  & -0.369431 & 1.179810   \\ 
		\hline
		4  & 2.610411 & 5.211440  & -0.328406 & 1.217261   \\ 
		\hline
		5  & 3.007662 & 6.006225  & -0.301045 & 1.142278   \\ 
		\hline
		6  & 3.380347 & 6.751789  & -0.281280 & 0.986919   \\ 
		\hline
		7  & 3.733587 & 7.458411  & -0.266220 & 0.770813   \\ 
		\hline
		8  & 4.070887 & 8.133123  & -0.254297 & 0.506915   \\ 
		\hline
		9  & 4.394782 & 8.781000  & -0.244582 & 0.204272   \\ 
		\hline
		10 & 4.707173 & 9.405855  & -0.236486 & -0.130527  \\ 
		\hline
		11 & 5.009538 & 10.010645 & -0.229617 & -0.492519  \\ 
		\hline
		12 & 5.303051 & 10.597720 & -0.223702 & -0.877859  \\ 
		\hline
		13 & 5.588663 & 11.168989 & -0.218546 & -1.283506  \\ 
		\hline
		14 & 5.867162 & 11.726026 & -0.214003 & -1.707004  \\ 
		\hline
		15 & 6.139207 & 12.270150 & -0.209966 & -2.146341  \\ 
		\hline
		16 & 6.405358 & 12.802481 & -0.206349 & -2.599845  \\ 
		\hline
		17 & 6.666093 & 13.323979 & -0.203088 & -3.066108  \\ 
		\hline
		18 & 6.921829 & 13.835474 & -0.200128 & -3.543935  \\ 
		\hline
		19 & 7.172927 & 14.337693 & -0.197428 & -4.032299  \\ 
		\hline
		20 & 7.419707 & 14.831273 & -0.194953 & -4.530312  \\
		\hline
	\end{tabular}
	\caption{Generating orbits for fixed~$\theta= \pi /18$.}
	\label{table generating theta10}
\end{table}

\begin{table}[tb]
	\centering
	\begin{tabular}{|c||c|c|c|c|c|c|} 
		\hline
		& \multicolumn{2}{c|}{Sun-Jupiter} & \multicolumn{2}{c|}{Earth-Moon} & \multicolumn{2}{c|}{Pluto-Charon}  \\ 
		\hline
		$I$  & $q_0$      & $\mathcal{A}$               & $q_0$      & $\mathcal{A}$              & $q_0$      & $\mathcal{A}$                 \\ 
		\hline\hline
		1  & 2.288901  & 0.485253             & 2.281639  & 0.590492            & 2.216384  & 1.209834               \\ 
		\hline
		2  & 3.396803  & 0.984282             & 3.392470  & 1.089921            & 3.353351  & 1.733091               \\ 
		\hline
		3  & 4.350951  & 1.192265             & 4.347446  & 1.297962            & 4.315925  & 1.949195               \\ 
		\hline
		4  & 5.211184  & 1.229658             & 5.208155  & 1.335366            & 5.180883  & 1.990792               \\ 
		\hline
		5  & 6.005993  & 1.154636             & 6.003252  & 1.260344            & 5.978585  & 1.918388               \\ 
		\hline
		6  & 6.751575  & 0.999251             & 6.749047  & 1.104955            & 6.726282  & 1.764807               \\ 
		\hline
		7  & 7.458211  & 0.783125             & 7.455838  & 0.888823            & 7.434486  & 1.550010               \\ 
		\hline
		8  & 8.132933  & 0.519211             & 8.130687  & 0.624905            & 8.110475  & 1.287122               \\ 
		\hline
		9  & 8.780818  & 0.216555             & 8.778675  & 0.322245            & 8.759385  & 0.985285               \\ 
		\hline
		10 & 9.405681  & -0.118255            & 9.403624  & -0.012570           & 9.385118  & 0.651146               \\ 
		\hline
		11 & 10.010477 & -0.480255            & 10.008494 & -0.374574           & 9.990653  & 0.289708               \\ 
		\hline
		12 & 10.597558 & -0.865603            & 10.595640 & -0.759926           & 10.578382 & -0.095163              \\ 
		\hline
		13 & 11.168832 & -1.271256            & 11.166970 & -1.165583           & 11.150223 & -0.500404              \\ 
		\hline
		14 & 11.725873 & -1.694760            & 11.724061 & -1.589089           & 11.707773 & -0.923547              \\ 
		\hline
		15 & 12.270001 & -2.134103            & 12.268235 & -2.028434           & 12.252356 & -1.362572              \\ 
		\hline
		16 & 12.802335 & -2.587611            & 12.800611 & -2.481946           & 12.785104 & -1.815798              \\ 
		\hline
		17 & 13.323836 & -3.053878            & 13.322149 & -2.948216           & 13.306981 & -2.281813              \\ 
		\hline
		18 & 13.835334 & -3.531709            & 13.833682 & -3.426048           & 13.818826 & -2.759415              \\ 
		\hline
		19 & 14.337556 & -4.020076            & 14.335935 & -3.914418           & 14.321365 & -3.247575              \\ 
		\hline
		20 & 14.831139 & -4.518092            & 14.829547 & -4.412436           & 14.815244 & -3.745402              \\
		\hline
	\end{tabular}
\caption{Continued orbits for fixed~$\theta = \pi/18$ for astronomical mass ratios.}
\label{table continued astro theta10}
\end{table}

\begin{table}[tb]
	\centering
	\begin{tabular}{|c||c|c|c|c|} 
		\hline
		& \multicolumn{2}{c|}{$\mu=0.2$} & \multicolumn{2}{c|}{$\mu=0.5$}  \\ 
		\hline
		$I$  & $q_0$      & $\mathcal{A}$          & $q_0$      & $\mathcal{A}$           \\ 
		\hline\hline
		1  & 2.150085  & 1.693790        & 1.816536  & 3.069224         \\ 
		\hline
		2  & 3.313535  & 2.246562        & 3.137846  & 3.758241         \\ 
		\hline
		3  & 4.283998  & 2.472591        & 4.146060  & 4.022613         \\ 
		\hline
		4  & 5.153181  & 2.519383        & 5.033830  & 4.088841         \\ 
		\hline
		5  & 5.953545  & 2.450237        & 5.846074  & 4.031746         \\ 
		\hline
		6  & 6.703154  & 2.298919        & 6.603994  & 3.888747         \\ 
		\hline
		7  & 7.412800  & 2.085796        & 7.319979  & 3.681770         \\ 
		\hline
		8  & 8.089940  & 1.824206        & 8.002119  & 3.424935         \\ 
		\hline
		9  & 8.739790  & 1.523409        & 8.656086  & 3.127945         \\ 
		\hline
		10 & 9.366317  & 1.190125        & 9.286061  & 2.797789         \\ 
		\hline
		11 & 9.972531  & 0.829403        & 9.895240  & 2.439691         \\ 
		\hline
		12 & 10.560854 & 0.445144        & 10.486136 & 2.057670         \\ 
		\hline
		13 & 11.133216 & 0.040432        & 11.060773 & 1.654893         \\ 
		\hline
		14 & 11.691232 & -0.382250       & 11.620811 & 1.233905         \\ 
		\hline
		15 & 12.236233 & -0.820866       & 12.167634 & 0.796784         \\ 
		\hline
		16 & 12.769361 & -1.273728       & 12.702410 & 0.345255         \\ 
		\hline
		17 & 13.291585 & -1.739416       & 13.226138 & -0.119236        \\ 
		\hline
		18 & 13.803747 & -2.216723       & 13.739680 & -0.595463        \\ 
		\hline
		19 & 14.306580 & -2.704616       & 14.243785 & -1.082372        \\ 
		\hline
		20 & 14.800730 & -3.202198       & 14.739113 & -1.579057        \\
		\hline
	\end{tabular}
\caption{Continued orbits for fixed~$\theta = \pi/18$ for nonastronomical mass ratios.}
\label{table continued nonastro theta10}
\end{table}

\clearpage